\documentclass[a4paper,10pt]{article}
\usepackage[english]{babel}
\usepackage{a4wide}
\usepackage{enumitem}
\usepackage{multicol}
\usepackage[utf8]{inputenc}

\usepackage{amsmath,amssymb,amsthm,array,textpos,float,dsfont,verbatim,fancyvrb}
\usepackage{breqn}
\usepackage[dvipsnames]{xcolor}
\usepackage{mathtools,makecell,graphicx,wrapfig,caption,subcaption,tikz,tikz-cd,framed}
\tikzset{
  symbol/.style={
    draw=none,
    every to/.append style={
      edge node={node [sloped, allow upside down, auto=false]{$#1$}}}
  }
}
\PassOptionsToPackage{hyphens}{url}\usepackage{hyperref}

\theoremstyle{definition}
\newtheorem{definition}{Definition}[section]
\newtheorem{example}[definition]{Example}

\theoremstyle{plain}
\newtheorem{theorem}[definition]{Theorem}
\newtheorem{proposition}[definition]{Proposition}
\newtheorem{corollary}[definition]{Corollary}
\newtheorem{lemma}[definition]{Lemma}

\theoremstyle{remark}
\newtheorem{remark}[definition]{Remark}

\newtheorem{notation}[definition]{Notation}

\numberwithin{equation}{section}

\newcommand\qt[1]{\quad \text{#1} \quad}
\newcommand\st[1]{\text{ #1 }}

\newcommand\Aut[1]{\text{Aut}\left(#1\right)}
\newcommand\End[1]{\text{End}\left(#1\right)}
\newcommand\Ker[1]{\text{Ker}\left(#1\right)}

\newcommand\Image[1]{\text{Im}\left(#1\right)}

\newcommand\Fix[1]{\text{Fix}(#1)}
\newcommand\MF[1]{\text{MF}(#1)}

\newcommand\Spec[1]{\text{Spec\textsubscript{R}}\left(#1\right)}

\makeatletter
\newcommand{\biggg}{\bBigg@\thr@@}
\newcommand{\Biggg}{\bBigg@{3}}

\makeatother

\makeatletter
\renewcommand*\env@matrix[1][*\c@MaxMatrixCols c]{%
  \hskip -\arraycolsep
  \let\@ifnextchar\new@ifnextchar
  \array{#1}}
\makeatother

\title{The Reidemeister spectrum of 2-step nilpotent groups determined by graphs}
\author{Karel Dekimpe\thanks{Research supported by long term structural funding - Methusalem grant of the Flemish Government.}, Maarten Lathouwers\thanks{Researcher funded by FWO PhD-fellowship fundamental research (file number: 1102422N).}}
\date{\today}

\begin{document}
\maketitle
\begin{abstract}
In this paper we study the Reidemeister spectrum of 2-step nilpotent groups associated to graphs. We develop three methods, based on the structure of the graph, that can be used to determine the Reidemeister spectrum of the associated group in terms of the Reidemeister spectra of groups associated to smaller graphs. We illustrate our methods for several families of graphs, including all the groups associated to a graph with at most four vertices. We also apply our results in the context of topological fixed point theory for nilmanifolds. 
\end{abstract}


\section{Introduction}
In this paper we will be studying Reidemeister numbers (this is the number of so called twisted conjugacy classes) of automorphisms of a wide class of 2-step nilpotent groups. Twisted conjugacy finds its origin in topological fixed point theory (see below) but pops up in several branches of mathematics, such as representation theory (\cite{OnishchikArkadijL1990Lgaa}, \cite{SpringerT.A2006Tcis}), Galois cohomology (\cite{SerreJean-Pierre2002Gc}), cryptography (\cite{GryakJonathan2016Tsop}), $\dots$

Our own main motivation comes from topological fixed point theory, more specifically from Reidemeis\-ter-Nielsen fixed point theory. We give a short overview of the main aspects of this theory and refer the reader to \cite{jm06-1,jiang1983lectures,Kiang} for more details.

Let $f:X\to X$ be a map on a closed manifold $X$ and denote with $\Fix{f}=\{x \in X \mid f(x) =x\}$ the set of fixed points of $f$. The main objective of Reidemeister-Nielsen fixed point theory is to find a good estimate for the minimal value of $\#\Fix{g}$ where $g$ is a map which is homotopic to $f$. Let us call this value $\MF{f}$.

To study the fixed points of $f$ one considers the universal covering space $p: \tilde{X} \to X$ of $X$. Then $f$ can be lifted to a map $\tilde{f}: \tilde{X} \to \tilde{X}$ (with $p \circ \tilde{f}= f \circ p$) and it is easy to see that $p(\Fix{\tilde{f}})\subseteq \Fix{f}.$
In fact, $\Fix{f}$ is the union of all $p(\Fix{\tilde{f}})$ where the union is taken over all possible lifts $\tilde{f}$ of $f$. The group of covering transformations of the universal covering is isomorphic to the fundamental group of $X$ and so we denote the group of covering transformations by $\pi(X)$. For any $\alpha, \beta \in \pi(X)$ and any lift $\tilde{f}$ of $f$ it holds that $\alpha \circ \tilde{f} \circ \beta$ is again a lift of $f$. It follows that $\pi(X)$ acts on the set of all lifts $\tilde{f}$ of $f$ via conjugation, so $\gamma\cdot \tilde{f} = \gamma \circ \tilde{f} \circ \gamma^{-1}$. We denote the orbit of $\tilde{f}$ by $[\tilde{f}]$ and call this the lifting class of $\tilde{f}$. Then we have that 
for all $\tilde{f}'\in [\tilde{f}]$ it holds that $p(\Fix{\tilde{f}'})= p(\Fix{\tilde{f}})$ while   
$p(\Fix{\tilde{f}'})\cap p(\Fix{\tilde{f}})=\emptyset$ in case $[\tilde{f'}]\neq [\tilde{f}]$. From this we can conclude that 
\[ \Fix{f} = \bigcup_{[\tilde{f}]} p(\Fix{\tilde{f}}),\]
which is a disjoint union. So in this union we consider one subset $p(\Fix{\tilde{f}})$ for each lifting class $[\tilde{f}]$ and we call this the fixed point class of $f$ determined by the lifting class $[\tilde{f}]$. Note that a fixed point class can be empty, but we still consider two empty fixed point classes different in case they are determined by a different lifting class. Hence, the number of fixed point classes is the same as the number of lifting classes and this number is called the Reidemeister number of $f$ and is denoted by $R(f)$. 

There is an algebraic way to count the fixed point classes of a map $f$ and this goes as follows. Fix one lifting $\tilde{f}_0$ of $f$. Then any other lift of $f$ can be written uniquely as a composition $\alpha \circ \tilde{f}_0$ for some $\alpha \in \pi(X)$. So the set of liftings of $f$ is in one-to-one correspondence with the fundamental group $\pi(X)$. The lift $\tilde{f}_0$ determines an endomorphism $f_\ast$ of $\pi(X)$ by the relation $f_\ast(\alpha) \circ \tilde{f}_0 = \tilde{f}_0 \circ \alpha$ (for all $\alpha\in \pi(X))$. Note that under the right identification of $\pi(X)$ with the fundamental group of $X$, $f_\ast$ is just the usual induced endomorphism of $f$ on the fundamental group of $X$.

Now $\alpha \circ \tilde{f}_0= \gamma \circ (\beta \circ \tilde{f}_0) \circ \gamma^{-1}$ if and only if 
$\alpha= \gamma \circ\beta \circ f_\ast(\gamma^{-1})$. In this case, we will say that $\alpha$ and $\beta$ are twisted conjugate with respect to $f_\ast$. Being twisted conjugate is an equivalence relation on $\pi(X)$ and the number of equivalence classes is called the Reidemeister number of the morphism $f_\ast$ and is denoted by $R(f_\ast)$. From the above we have that $R(f)=R(f_\ast)$ and so counting twisted conjugacy classes is an algebraic way of counting fixed point classes (or lifting classes). 

Although the Reidemeister number of a map gives already some information about the fixed point classes of $f$, in general this number does not really give information on $\MF{f}$, the minimal number of fixed points in the homotopy class of $f$. 
There is a second number, the Nielsen number of $f$ which does provide more information, but unfortunately is much more difficult to compute in general. To define the Nielsen number, there is a way to attach to each fixed point class an index, which is an integer. It would lead us to far to explain this index in more detail, but the idea is that a fixed point class has index 0 if it can disappear (become empty) via a homotopy. A fixed point class is called essential (resp.\ non essential) if it has an index $\neq 0$ (resp.\ $=0$). The Nielsen number of $f$, denoted by $N(f)$, is then the number of essential fixed point classes of $f$. This Nielsen number (and also the Reidemeister number) is a homotopy invariant and by a result of Wecken (\cite{wecken1942fixpunktklassen}) it is known that $MF(f)=N(f)$ for all manifolds of dimension at least 3. 

The focus of this paper lies on the class of nilmanifolds, these are obtained as quotient spaces $X=N\backslash G$, where $G$ is a simply connected nilpotent Lie group and $N$ is a uniform lattice of $G$. Such a uniform lattice $N$ is a finitely generated torsion-free nilpotent group and $\pi(X)=N$ completely determines the nilmanifold $X=N \backslash G$ up to diffeomorphism (see e.g.\ \cite{raghunathan1972discrete,ov93-1}).  

For the class of nilmanifolds, there is a very strong relation between the Reidemeister number and the Nielsen number of a map $f$ on such a manifold. Indeed, we have that (see \cite{hk97-1}):
\[ \left\{\begin{array}{lcl}
N(f) = R(f) & \Longleftrightarrow & R(f)<\infty \\
N(f) = 0 & \Longleftrightarrow & R(f) = \infty
\end{array}\right. .\]
As a conclusion we see that for nilmanifolds, we obtain a full understanding of the minimal number of fixed points ($\MF{f}=N(f)$) in the homotopy class of a map $f$ by studying the Reidemeister number $R(f)$ of that map and hence by studying the number $R(f_\ast)$ of twisted conjugacy classes of the induced endomorphism $f_\ast$.

\medskip

It is not so difficult to see that for any nilmanifold $N\backslash G$ and any non negative integer $n$ there is a self map $f$ of $N\backslash G$ with 
$N(f)=n$ (\cite[Theorem 6.1]{dekimpe2020}). The situation for self homeomorphisms (or self homotopy equivalences) is much more subtle and corresponds to the case where $f_\ast$ is an automorphism of $N$. The set of all possible Reidemeister numbers one can obtain for these self homotopy equivalences $f$ (automorphisms $f_\ast$) is called the Reidemeister spectrum of the manifold (or of the fundamental group $N$).
In this paper we study this Reidemeister spectrum for groups $N$ which are 2-step nilpotent and are associated to a graph. The nilmanifolds with such a 2-step nilpotent fundamental group have been the object of study in many geometric contexts and form a rich family of interesting examples. We refer to \cite{voorbeeld0, voorbeeld1, voorbeeld2, voorbeeld3} for some recent examples in this direction. Moreover, results on 2-step nilpotent groups can be used to study general nilpotent groups by considering their 2-step nilpotent quotient.

In the next section we recall some preliminaries on nilpotent groups and twisted conjugacy.
Thereafter, we describe the class of 2-step nilpotent groups associated to a graph. In the next three sections we develop each time a general method, based on the structure of the graph, that can be used to determine the Reidemeister spectrum of the associated group, by reducing it to the situation of smaller graphs. We end by illustrating our methods in some general examples and give a full list of Reidemeister spectra for all graphs with at most 4 vertices.



\section{Preliminaries on nilpotent groups and twisted conjugacy}
\subsection{Nilpotent groups}
For any group $G$ we denote with $\gamma_i(G)$ (for $i\in \mathbb{N}_0$) the \textit{lower central series of $G$}, i.e. the nested series of subgroups of $G$ defined by $\gamma_1(G):=G$ and $\gamma_{i+1}(G):=[\gamma_i(G),G]$ (for $i\in \mathbb{N}_0$). The group $G$ is said to be \textit{$c$-step nilpotent} if $\gamma_c(G)\neq 1$ and $\gamma_{c+1}(G)=1$. It is generally known (see e.g. \cite[Theorem 17.2.2]{kargapolov1979}) that any finitely generated nilpotent group $G$ has a series $1=G_1\lhd G_2\lhd \dots \lhd G_s=G$ with cyclic factors, i.e. $G_{i+1}/G_i$ is cyclic for any $i=1,2,\dots,s-1$. The \textit{Hirsch number} $h(G)$ of $G$ is the number of infinite cyclic factors in such a series. The next lemma describes some properties of the Hirsch number.

\begin{lemma}[{\cite[page 16]{segal1983polycyclic}}]\label{lemma:HN properties}
If $G$ is a finitely generated nilpotent group, then the Hirsch number is well-defined (meaning that it is independent of the choice of series of $G$ with cyclic factors). If $H\subseteq G$ is a subgroup of $G$ and $N\lhd G$ a normal subgroup, then the following holds:
\begin{enumerate}
\item[(i)] $h(H)\leq h(G)$
\item[(ii)] $h(H)=h(G)\: \Longleftrightarrow \: [G:H]<\infty$
\item[(iii)] $h(G)=h(N)+h(G/N)$
\item[(iv)] $h(G)=0\: \Longleftrightarrow\: |G|\:<\infty$
\end{enumerate}
\end{lemma}

In section \ref{sec:associating groups to graphs} we describe how to associate a finitely generated torsion-free $2$-step nilpotent group to any finite undirected simple graph. The next two lemmas will be used frequently.

\begin{lemma}[\cite{KarrassAbraham1976Cgt:}]\label{lemma:bilinearity [.,.] 2-nilpotent}
If $G$ is a $2$-step nilpotent group, then $[\:.\:,.\:]:G\times G\to G$ is bilinear, i.e. for all $g_1,g_2,g_1',g_2'\in G$ it holds that 
\[[g_1g_2,g_1'g_2']=[g_1,g_2']\cdot[g_2,g_2']\cdot[g_1,g_1']\cdot[g_2,g_1'].\]
\end{lemma}

\begin{lemma}\label{lemma:surj ==> aut}
If $G$ is a finitely generated torsion-free nilpotent group and $\varphi:G\to G$ is a surjective morphism, then $\varphi$ is also injective. In particular, $\varphi$ is an automorphism of $G$.
\end{lemma}
\begin{proof}
Since $\varphi$ is surjective, it follows by the first isomorphism theorem and Lemma \ref{lemma:HN properties} (iii) that
\[h(\Ker{\varphi})=h(G)-h(G/\Ker{\varphi})=h(G)-h(\Image{\varphi})=0.\]
So Lemma \ref{lemma:HN properties} (iv) implies that $\Ker{\varphi}$ is finite. However, since $G$ is torsion-free we obtain that $\Ker{\varphi}$ is a finite torsion-free group. Hence, $\Ker{\varphi}$ is trivial.
\end{proof}

Fix any endomorphism $\varphi\in \End{G}$. We denote with $\varphi_i\in \End{\gamma_i(G)/\gamma_{i+1}(G)}$ (for $i\in\mathbb{N}_0$) the induced morphisms on the factors of the lower central series of $G$. Since the terms of the lower central series are characteristic subgroups of $G$, these induced morphisms are well-defined. Moreover, if $G$ is a finitely generated $c$-step nilpotent group, then these factors are finitely generated abelian groups (see for example Lemma 17.2.1 in \cite{kargapolov1979}). 

In case $\gamma_i(G)/\gamma_{i+1}(G)$ is a free abelian group, so isomorphic to $\mathbb{Z}^k$ for some $k$, we can describe 
$\varphi_i$ by a $k\times k$ matrix over $\mathbb{Z}$ and in this way we can talk about the eigenvalues of $\varphi_i$ and the determinant of 
$\varphi_i$. In case $\gamma_i(G)/\gamma_{i+1}(G)$ has torsion, we can still talk about the eigenvalues of $\varphi_i$ by which we mean the eigenvalues of the induced endomorphism on the torsion-free quotient $\frac{\gamma_i(G)/\gamma_{i+1}(G)}{\tau(\gamma_i(G)/\gamma_{i+1}(G))}$ (where $\tau(\gamma_i(G)/\gamma_{i+1}(G))$ denotes the torsion subgroup of $\gamma_i(G)/\gamma_{i+1}(G)$).

The next result will be needed later on.

\begin{lemma}\label{lemma:phi_1 aut ==> phi aut}
Let $G$ be a finitely generated torsion-free nilpotent group and $\varphi:G\to G$ a morphism such that $\varphi_1:G/\gamma_2(G)\to G/\gamma_2(G)$ is an automorphism. Then $ \varphi$ is an automorphism.
\end{lemma}

\begin{proof}
Since $\varphi_1$ is surjective, it follows that
\[\frac{G}{\gamma_2(G)}=\Image{\varphi_1}=\frac{\Image{\varphi}\gamma_2(G)}{\gamma_2(G)}.\]
So it holds that $\Image{\varphi}\gamma_2(G)=G$. Since $G$ is nilpotent, this implies tat $\Image{\varphi}=G$ (see e.g. \cite[Theorem 16.2.5]{kargapolov1979}). Hence, $\varphi$ is a surjective morphism and thus Lemma \ref{lemma:surj ==> aut} yields the result.
\end{proof}

\subsection{Twisted conjugacy}
In the introduction we already considered the notion of twisted conjugacy. Let us now define this concept in some more detail.
Let $G$ be a group and $\varphi\in \End{G}$ a group endomorphism. Two elements $a,b\in G$ are called \textit{$\varphi$-conjugate} or \textit{twisted conjugate} (denoted by $a\sim_{\varphi} b$) if there exists some third element $c\in G$ such that $a=cb\varphi(c)^{-1}$. This induces an equivalence relation $\sim_{\varphi}$ on $G$. The equivalence classes are called the \textit{Reidemeister classes of $\varphi$} and the number of equivalence classes is called the \textit{Reidemeister number $R(\varphi)$ of $\varphi$}. The collection of all Reidemeister numbers, by only considering automorphisms of $G$, is called the \textit{Reidemeister spectrum of $G$} and is denoted by $\Spec{G}$. Formally, we define $\Spec{G}:=\{R(\varphi)\:|\: \varphi\in \Aut{G}\}\subseteq \mathbb{N}_0\cup \{\infty\}$. If $\Spec{G}=\{\infty\}$, then $G$ is said to have the \textit{$R_{\infty}$--property}. If $\Spec{G}=\mathbb{N}_0\cup \{\infty\}$, then $G$ has a \textit{full Reidemeister spectrum}.

The Reidemeister spectrum of a group is in general difficult to compute. However, when working with finitely generated nilpotent groups, there are some well-known techniques to study their Reidemeister spectrum. 

We list two theorems without a proof.

\begin{proposition}[{\cite[Lemma 2.2]{DekimpeKarel2014TRip}} and {\cite[Corollary 4.2]{roman2011twisted}}]\label{prop:equivalent statements R-inf}
Let $G$ be a finitely generated $c$-step nilpotent group and $\varphi$ an automorphism of $G$. Then the following are equivalent:
\begin{enumerate}
\item[(i)] $R(\varphi)=\infty$
\item[(ii)] There exists some $i=1,2,\dots,c$ such that $\varphi_i$ has $1$ as an eigenvalue.
\item[(iii)] There exists some $i=1,2,\dots,c$ such that $R(\varphi_i)=\infty$.
\end{enumerate}
\end{proposition}

\begin{theorem}[{\cite[Proposition 5]{dekimpe2021r_}} and {\cite[Lemma 2.7]{roman2011twisted}}]\label{thm:RN product}
Let $G$ be a finitely generated nilpotent group. Let
\[G=G_1\supseteq G_2\supseteq G_3 \supseteq \dots \supseteq G_{c}\supseteq G_{c+1}=1\]
be a central series of $G$ and $\varphi\in \Aut{G}$ such that the following holds:
\begin{enumerate}
\item All the factors $G_i/G_{i+1}$ (with $i=1,2,\dots,c$) are torsion-free.
\item For all terms $G_i$ (with $i=1,2,\dots,c+1$) it holds that $\varphi(G_i)=G_i$.
\end{enumerate}
Then it holds that:
\[R(\varphi)=\prod_{i=1}^c R(\varphi_i)\]
where $\varphi_i:G_i/G_{i+1}\to G_i/G_{i+1}$ (with $i=1,2,\dots,c$) are the induced automorphisms on the factor groups $G_i/G_{i+1}$.
\end{theorem}

\begin{remark}
Note that in Theorem \ref{thm:RN product} we used the notation $\varphi_i$ to denote the induced automorphisms on the factor groups of the given central series of $G$. However, whenever we do not mention a central series, we reserve this notation for the induced automorphisms on the factors of the lower central series (as we introduced above Lemma \ref{lemma:phi_1 aut ==> phi aut}).
\end{remark}

When describing the Reidemeister spectrum of a group we will frequently use the map $|\:\cdot\:|_{\infty}$ which is defined by
\[|\:\cdot\:|_{\infty}:\mathbb{Z}\to \mathbb{N}_0\cup \{\infty\}:x\mapsto |x|_{\infty}:=\begin{cases}
|x| &\st{if} x\neq 0\\
\infty &\st{if} x=0 \end{cases}\]
where $|\:\cdot\:|$ denotes the absolute value.


\section{2-step nilpotent groups associated to graphs}\label{sec:associating groups to graphs}
To any undirected finite simple graph $\Gamma(V,E)$, we can associate a finitely generated 2-step nilpotent group. We do this by considering the $2$-step nilpotent quotient of the \textit{right angled Artin group} associated to $\Gamma$. More precisely, the idea is that we take the vertices (which we denote with $x_i$) as generators of our group and require that two of these generators commute if the corresponding vertices are connected via an edge. We denote with $y_{i,j}=[x_j,x_i]$ the commutators of two vertices that are not connected via an edge. At last, we make the group $2$-step nilpotent by adding the constraints that the $y_{i,j}$ commute with all the vertices. The formal definition is given below.

\begin{definition}\label{def:group associated to graph}
Let $\Gamma(\{x_1,x_2,\dots,x_n\},E)$ be an undirected finite simple graph. We define the group $G_\Gamma$ by setting
\[G_{\Gamma}=\Bigg \langle \begin{array}{ll} x_1,x_2,\dots,x_n,\\
y_{i,j} \st{if} x_ix_j\not\in E \st{and} i<j \end{array}\: \Biggg\vert \:   
\begin{array}{ll}
[x_j,x_i]=1 &\st{if} x_ix_j\in E \\
\relax [x_j,x_i]=y_{i,j} &\st{if} x_ix_j\not\in E \st{and} i<j \\
\relax [x_l,y_{i,j}]=1 &\: l=1,2,\dots, n;\: x_ix_j\not\in E \st{and} i<j
\end{array}
\Bigg \rangle\]
\end{definition}

From now on, we will use $\Gamma$ to denote an undirected finite simple graph. Note that if $\Gamma$ is the complete graph on $n\in \mathbb{N}_0$ vertices, then the associated group $G_{\Gamma}$ is isomorphic with $\mathbb{Z}^n$. The Reidemeister spectrum of $\mathbb{Z}^n$ is well-known (see for example \cite[Section 3]{roman2011twisted}) and is given by
\begin{equation}\label{eq:Spec Z^n}
\Spec{\mathbb{Z}^n}=
\begin{cases}
\{2,\infty\} &\st{if} n=1\\
\mathbb{N}_0\cup \{\infty\} &\st{if} n\geq 2
\end{cases}.
\end{equation}
If on the other hand $\Gamma$ is the graph on $n$ vertices without any edges, then the associated group $G_{\Gamma}$ is isomorphic with $F_n/\gamma_3(F_n)$ (where $F_n$ denotes the free group of rank $n$). These groups are also known as the \textit{free nilpotent groups of rank $n$ and nilpotency class 2} and are frequently denoted by $N_{n,2}$. K.~Dekimpe, S.~Tertooy and A.R.~Vargas extended in \cite[Section 4]{dekimpe2020} the result from V.~Roman'kov (see \cite[Section 3]{roman2011twisted}) to
\begin{equation}\label{eq:Spec N_n,2}
\Spec{N_{n,2}}=\begin{cases}
2\mathbb{N}_0\cup \{\infty\} &\text{if } n=2\\
(2\mathbb{N}_0-1)\cup 4\mathbb{N}_0\cup \{\infty\}&\text{if } n=3\\
\mathbb{N}_0\cup \{\infty\} &\text{if } n\geq 4
\end{cases}.
\end{equation}

Denote with $N:=|\{(i,j)\:|\: x_ix_j\not\in E \st{and} i<j\}|$. To simplify notation, we fix an order $y_1,y_2,\dots, y_N$ to denote the elements $y_{i,j}$ where we define $y_l:= y_{i_l,j_l}$. Using the definition, it follows that any element of $G_{\Gamma}$ can be uniquely written as $x_1^{z_1}x_2^{z_2}\dots x_n^{z_n}y_1^{t_1}y_2^{t_2}\dots y_N^{t_N}$ with $z_i, t_l\in \mathbb{Z}$ and that the multiplication in $G_{\Gamma}$ is given by
\[(x_1^{z_1}\dots x_n^{z_n}y_1^{t_1}\dots y_N^{t_N})(x_1^{v_1}\dots x_n^{v_n}y_1^{s_1}\dots y_N^{s_N})=x_1^{z_1+v_1}\dots x_n^{z_n+v_n}y_1^{t_1+s_1+v_{i_1}z_{j_1}}\dots y_N^{t_N+s_N+v_{i_N}z_{j_N}}\]
for $z_i,v_i,t_l,s_l\in \mathbb{Z}$.

Using the operation in $G_{\Gamma}$, we obtain expressions for the center and commutator subgroup of $G_{\Gamma}$.

\begin{lemma} \label{lemma:factor groups Ggamma}
With the notations from above, we have that:
\begin{align*}
Z(G_\Gamma)&= \bigtimes_{i=1}^N \left<y_i\right>\times \bigtimes_{\substack{i=1,2,\dots,n;\\ \deg(x_i)=n-1}} \left<x_i\right>\cong \mathbb{Z}^{N+ (\text{the number of vertices of degree } n-1)}\\
\st{and} \gamma_2(G_\Gamma)&=\bigtimes_{i=1}^N \left<y_i\right>
\cong \mathbb{Z}^N
\end{align*}
\end{lemma}

Note that Lemma \ref{lemma:factor groups Ggamma} implies that if $\Gamma$ is not a complete graph, then the associated group $G_{\Gamma}$ is a finitely generated torsion-free 2-step nilpotent group. Hence, we can apply Theorem \ref{thm:RN product} to the lower central series of $G_{\Gamma}$. Combined with the well-known description of the Reidemeister spectrum for finitely generated torsion-free abelian groups (see e.g. \cite{gonccalves2009twisted}) we obtain the next result.

\begin{lemma} \label{lemma:Ggamma RN product induced aut}
For any $\varphi\in \Aut{G_\Gamma}$ we have that:
\[R(\varphi)=R(\varphi_1)R(\varphi_2)=|\det(Id-\varphi_1)|_{\infty}\:|\det(Id-\varphi_2)|_{\infty}\]
where $\varphi_1:G_\Gamma/\gamma_2(G_\Gamma)\to G_\Gamma/\gamma_2(G_\Gamma)$ and $\varphi_2:\gamma_2(G_\Gamma)\to \gamma_2(G_\Gamma)$ are the induced automorphisms on the factors of the lower central series.
\end{lemma}
We abuse notation and also denote with $\varphi_1$ (respectively $\varphi_2$) the matrix corresponding to the map $\varphi_1$ (respectively $\varphi_2$).

It is clear that if two graphs are isomorphic, then the associated groups are isomorphic. Also the converse is true. This can be proven by using the argument for the associated \textit{right-angled Artin groups} in \cite{droms1987isomorphisms}. In his argument Droms actually shows that when the 2-step nilpotent quotients of the right-angled Artin groups associated to the two graphs are isomorphic, that it follows that the graphs are isomorphic. However, since $G_{\Gamma}$ is precisely the 2-step nilpotent quotient of the right-angled Artin group associated to $\Gamma$, this argument suffices to conclude the proof of the following Lemma.

\begin{lemma}\label{lemma:isomorphic graphs}
Let $\Gamma_1$ and $\Gamma_2$ be two undirected simple graphs. The graphs $\Gamma_1$ and $\Gamma_2$ are isomorphic if and only if $G_{\Gamma_1}\cong G_{\Gamma_2}$.
\end{lemma}


\section{The first method: degree of the vertices}
Recall that we want to determine the Reidemeister spectrum of the finitely generated 2-step nilpotent groups associated to graphs. To do so we develop three methods. For the first method we describe characteristic subgroups based on the degree of the vertices. The other two methods allow us to partition the graph by using the \textit{simplicial join} or the \textit{disjoint union}.

\medskip

For any graph $\Gamma(\{x_1,x_2,\dots,x_n\},E)$ and any $d\in \{1,2,\dots,n-1\}$ we define the vertex set
\[V_d:= \{x_i\: |\: i=1,2,\dots,n \st{and} \deg{x_i}\geq d\}\]
and the subgroup $H_d\subseteq G_{\Gamma}$ by
\[H_d:= \left \{\prod_{x_i\in V_d} x_i^{z_i}\: \prod_{l=1}^N y_l^{t_l}\: \middle\vert \: z_i,t_l\in \mathbb{Z} \right \}.\]
We argue that these subgroups $H_d$ are characteristic subgroups of $G_{\Gamma}$. For this, we need two lemmas that describe the Hirsch number of the centralizers of elements of $G_{\Gamma}$.

\begin{lemma}\label{lemma:HN centralizer 1 vertex}
For any $z_{i_0},t_l\in \mathbb{Z}$ (for some $i_0=1,2,\dots,n$ and any $l=1,2,\dots, N$) with $z_{i_0}\neq 0$ it holds that
\[h\left(Z_{G_{\Gamma}}\left( x_{i_0}^{z_{i_0}}\prod_{l=1}^N y_l^{t_l}\right) \right )=\deg(x_{i_0})+N+1.\]
\end{lemma}

\begin{proof}
By the operation in $G_{\Gamma}$, it follows that the centralizer of $x:=x_{i_0}^{z_{i_0}}\prod_{l=1}^N y_l^{t_l}$ has the following form
\[Z_{G_{\Gamma}}(x)=\{g\in G_{\Gamma}\: |\: [g,x]=1\}=\left\{\prod_{i=1}^n x_i^{v_i}\: \prod_{l=1}^N y_l^{s_l}\: \middle\vert \: \begin{array}{l} v_i,s_l\in \mathbb{Z}\\z_{i_0}v_i=0 \st{if} x_{i_0}x_i\not\in E\st{and} i_0\neq i\end{array}\right\}.\]
Since $z_{i_0}\neq 0$, this precisely means that
\[Z_{G_{\Gamma}}(x)=\left\{\prod_{\substack{i=1,2,\dots,n;\\x_{i_0}x_i\in E\st{or} i=i_0}} \hspace{-15pt} x_i^{v_i}\: \prod_{l=1}^N y_l^{s_l}\: \middle\vert \: v_i,s_l\in \mathbb{Z}\right\}.\]
Using Lemma \ref{lemma:HN properties} it follows that
\[h(Z_{G_{\Gamma}}(x))=N+|\{i=1,2,\dots,n\: |\: x_{i_0}x_i\in E\st{or}i=i_0\}|=\deg(x_{i_0})+N+1.\]
\end{proof}

\begin{lemma}\label{lemma:HN centralizer general element}
For any $z_i,t_l\in \mathbb{Z}$ (for $i=1,2,\dots,n$ and $l=1,2,\dots, N$) it holds that
\[h\left(Z_{G_{\Gamma}}\left( \prod_{i=1}^n x_i^{z_i}\: \prod_{l=1}^N y_l^{t_l}\right) \right )\leq \min\{n-1,\min_{\substack{i=1,2,\dots,n;\\ z_i\neq 0}}\deg(x_i)\}+N+1.\]
\end{lemma}

\begin{proof}
Denote $x:=\prod_{i=1}^n x_i^{z_i}\: \prod_{l=1}^N y_l^{t_l}$. If all $z_i$ are equal to zero, then $Z_{G_{\Gamma}}(x)=G_{\Gamma}$ and thus $h(Z_{G_{\Gamma}})=n+N$ and the result follows. So suppose that not all $z_i$ are zero. We can assume without loss of generality that
\[z_i\neq 0 \quad \Longleftrightarrow\quad i\in \{1,2,\dots,k\}\]
for some $k\in \{1,2,\dots,n\}$ and that
\[\deg(x_1)=\min_{i=1,2,\dots,k}\deg(x_i).\]
We denote with $m$ the number of $x_j$'s (with $j=2,3,\dots,k$) such that $x_1x_j\in E$. Without loss of generality, we can assume that
\[x_1x_j\not\in E \st{(with} j\in \{2,3,\dots,k\}) \quad \Longleftrightarrow\quad j=2,3,\dots,k-m.\]
By the operation in $G_{\Gamma}$, we obtain that
\[Z_{G_{\Gamma}}(x)=\left\{\prod_{i=1}^n x_i^{v_i}\: \prod_{l=1}^N y_l^{s_l}\: \middle\vert \: 
\begin{array}{ll} 
v_i,s_l\in \mathbb{Z}& \\
z_1v_j=v_1z_j\: &(\forall j=2,3,\dots,k-m)\\
z_iv_j=v_iz_j\: &(\forall i,j=2,3,\dots,k\st{with} x_ix_j\not\in E\st{and} i\neq j)\\
z_iv_j=0\: &(\forall i=1,2,\dots,k;j=k+1,\dots,n \st{with} x_ix_j\not\in E)\\
\end{array}\right\}.\]
We define the subgroup $H\subseteq G_{\Gamma}$ (that contains $Z_{G_{\Gamma}}(x)$) by setting
\[H:=\left\{\prod_{i=1}^n x_i^{v_i}\: \prod_{l=1}^N y_l^{s_l}\: \middle\vert \: 
\begin{array}{ll} 
v_i,s_l\in \mathbb{Z}& \\
z_1v_j=v_1z_j\: &(\forall j=2,3,\dots,k-m)\\
z_iv_j=0\: &(\forall i=1,2,\dots,k;j=k+1,\dots,n \st{with} x_ix_j\not\in E)\\
\end{array}\right\}.\]
By Lemma \ref{lemma:HN properties} (i) it suffices to argue that $h(H)\leq \deg(x_1)+N+1$. Since $z_1\neq 0$, the equations $z_1v_j=v_1z_j$ (with $j=2,3,\dots,k-m$) can only be satisfied if 
\[(v_1,v_2,\dots,v_{k-m})=\lambda\left( \frac{z_1}{d},\frac{z_2}{d},\dots,\frac{z_{k-m}}{d}\right)\quad \st{for some} \lambda\in \mathbb{Z}\]
where $d:=\gcd(z_1,z_2,\dots,z_{k-m})$. Since all the $z_l\neq 0$ (for $l=1,2,\dots,k$), it follows that $v_j=0$ for all $j=k+1,k+2,\dots,n$ with $x_ix_j\not\in E$ for some $i=1,2,\dots,k$. Hence, we obtain that
\[H=\left\{\left (\prod_{i=1}^{k-m} x_i^{z_i/d}\right)^{v_0}\prod_{i=k-m+1}^{k} x_i^{v_i}\: \prod_{j\in J}x_j^{v_j}\:\prod_{l=1}^N y_l^{s_l}\: \middle\vert \:
v_0,v_i,v_j,s_l\in \mathbb{Z}
\right\}\]
where we defined
\[J:=\{j=k+1,k+2,\dots,n\: |\: x_ix_j\in E \st{for all} i=1,2,\dots,k\}.\]
By using Lemma \ref{lemma:HN properties} we can indeed conclude that
\begin{align*}
h(H)&=1+m+|J|+N\leq 1+m+N+|\{j=k+1,k+2,\dots,n\: |\: x_1x_j\in E\}|\\
&=1+m+N+|\{j=2,3,\dots,n\: |\: x_1x_j\in E\}|-m\\
&=\deg(x_1)+N+1.
\end{align*}
\end{proof}

\begin{theorem}\label{thm:degree vertices}
The subgroups $H_d$ (for any $d=1,2,\dots, n-1$) are characteristic subgroups of $G_{\Gamma}$.
\end{theorem}

\begin{proof}
Fix any $d=1,2,\dots,n-1$, an automorphism $\varphi\in \Aut{G_{\Gamma}}$ and any $x_{i_0}\in V_d$. Suppose by contradiction that $\varphi(x_{i_0})\not\in H_d$. Hence, there exists some $x_{i_1}\not\in V_d$ and $z_i,t_l\in \mathbb{Z}$ with $z_{i_1}\neq 0$ such that
\[\varphi(x_{i_0})=\prod_{i=1}^n x_i^{z_i}\: \prod_{l=1}^N y_l^{t_l}.\]
By using Lemma \ref{lemma:HN centralizer 1 vertex} and Lemma \ref{lemma:HN centralizer general element} we now obtain that
\begin{align*}
d+N+1&\leq \deg{x_{i_0}}+N+1=h(Z_{G_{\Gamma}}(x_{i_0}))=h(Z_{G_{\Gamma}}(\varphi(x_{i_0})))\leq \deg{x_{i_1}}+N+1<d+N+1
\end{align*}
which is a contradiction and thus $\varphi(x_{i_0})\in H_d$. By Lemma \ref{lemma:factor groups Ggamma} and since $\gamma_2(G_{\Gamma})$ is a characteristic subgroup of $G_{\Gamma}$, it now follows that $\varphi(H_d)\subseteq H_d$. We can use completely the same argument to argue that $\varphi^{-1}(H_d)\subseteq H_d$ and thus we can conclude that $\varphi(H_d)=H_d$.
\end{proof}

Theorem \ref{thm:degree vertices} provides extra information about the automorphisms of $G_{\Gamma}$. Therefore, we will use it frequently in the rest of the paper to determine the Reidemeister spectrum of groups associated to graphs and to develop new methods to do so. Moreover, Theorem \ref{thm:degree vertices} can be used to describe graphs for which the associated finitely generated 2-step nilpotent groups have the $R_{\infty}$--property.

\begin{theorem}\label{thm:graphs having R-infinity}
Let $\Gamma(\{x_1,x_2,\dots,x_n\},E)$ be an undirected simple graph which has maximal degree $n-2$ and for which this degree is attained only once, then $G_{\Gamma}$ has the $R_{\infty}$--property.
\end{theorem}

\begin{proof}
Take any $\varphi\in \Aut{G_{\Gamma}}$. Assume without loss of generality that $x_1$ is the one vertex having degree $n-2$ and that $x_1x_2\not\in E$. Since $H_{n-2}$ is a characteristic subgroup of $G_{\Gamma}$ (by Theorem \ref{thm:degree vertices}), it follows that $\varphi_1(x_1\gamma_2(G_{\Gamma}))=x_1^{\pm 1}\gamma_2(G_{\Gamma})$.

Fix any $j\in \{3,4,\dots,n\}$ and denote $\varphi(x_j)=\prod_{i=1}^n x_i^{z_i}\: \prod_{l=1}^N y_l^{t_l}$. Since $x_1x_j\in E$, we obtain (by Lemma \ref{lemma:bilinearity [.,.] 2-nilpotent}) that
\[1=\varphi([x_j,x_1])=\left [\prod_{i=1}^n x_i^{z_i},x_1^{\pm 1}\right ]=\prod_{i=1}^n[x_i,x_1]^{\pm z_i}=[x_2,x_1]^{\pm z_2}.\]
Hence, we obtain that $z_2=0$. Since this argument is valid for all $j\in \{3,4,\dots,n\}$ and since $\varphi$ is an automorphism, it follows that the matrix of $\varphi_1$ (with respect to $\{x_1\gamma_2(G_{\Gamma}),x_2\gamma_2(G_{\Gamma}),\dots,x_n\gamma_2(G_{\Gamma})\}$) has the following form:
\[\begin{pmatrix}[c|c|ccc]
\pm 1 & b_1 &  & C &  \\
\hline
0 & \pm 1 & 0 & \hdots & 0 \\ 
\hline
0 & b_3 &  &  &  \\ 
\vdots & \vdots &  & A &  \\ 
0 & b_n &  &  & 
\end{pmatrix}\]
where $B:=\begin{pmatrix}
b_1 & \pm 1 & b_3 & \hdots & b_n
\end{pmatrix}^\top\in \mathbb{Z}^{n\times 1}$, $A\in \mathbb{Z}^{(n-2)\times (n-2)}$ and $C\in \mathbb{Z}^{1\times (n-2)}$. Since $\varphi_1$ is an automorphism, it holds that $A\in \text{GL}_{n-2}(\mathbb{Z})$. Applying Lemma \ref{lemma:Ggamma RN product induced aut} yields
\[R(\varphi)=R(\varphi_1)R(\varphi_2)=|\pm 1-1|_{\infty}\: |\pm 1-1|_{\infty}\: |\det(A-\mathds{1}_{n-2})|_{\infty}\: R(\varphi_2).\]
If one of the two $\pm 1$ is equal to $1$, then we get that $R(\varphi)=\infty$ and the result follows. Hence, we can assume without loss of generality that the two $\pm 1$ in the matrix representation are both equal to $-1$. However, by Lemma \ref{lemma:bilinearity [.,.] 2-nilpotent} we now obtain that
\begin{align*}
\varphi_2([x_2,x_1])&=[x_1^{b_1}x_2^{-1}x_3^{b_3}\dots x_n^{b_n},x_1^{-1}]=
[x_1,x_1]^{-b_1}[x_2,x_1]^{(-1)(-1)}[x_3,x_1]^{-b_3}\dots[x_n,x_1]^{-b_n}\\
&=[x_2,x_1].
\end{align*}
This implies that $[x_2,x_1]$ is an eigenvector of $\varphi_2$ with corresponding eigenvalue $1$ and thus by Proposition \ref{prop:equivalent statements R-inf} it follows that $R(\varphi)=\infty$.
\end{proof}

\begin{example}\label{cor:4V first R_infty}
If $\Gamma$ is the graph in Figure $\ref{fig:Graph 4V first R_infty}$, then $G_{\Gamma}$ has the $R_{\infty}$--property.
\begin{figure}[H]
\centering
\begin{tikzpicture}
    \tikzstyle{vertex}=
      [circle,draw,minimum size=1.2em,inner sep=0.2];
      
    \node[vertex] (x1) at (0,1) {$x_1$};
    \node[vertex] (x2) at (1,1) {$x_2$};
    \node[vertex] (x3) at (1,0) {$x_3$};
    \node[vertex] (x4) at (0,0) {$x_4$};

    \draw (x1)--(x2)--(x3);
\end{tikzpicture}
\caption{Graph for which the associated group has the $R_{\infty}$--property.} \label{fig:Graph 4V first R_infty}
\end{figure}
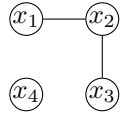
\end{example}

\begin{corollary}
Let $\Gamma(\{x_1,x_2,\dots,x_n\},E)$ be an undirected simple graph which has maximal degree $n-2$ and for which this degree is attained only once and let $M$ be the nilmanifold with fundamental group $G_\Gamma$. Then any self-homotopy equivalence of $M$ is homotopic to a fixed point free map. 
\end{corollary}

\begin{remark}\label{manifold dim 7}
The nilmanifold $M$ whose fundamental group is the group $G_\Gamma$ of Example \ref{cor:4V first R_infty} is an 8-dimensional (since $h(G_\Gamma)=8$) example of a 2-step nilmanifold for which every self-homotopy equivalence is homotopic to a fixed point free map. In Remark \ref{rem:cycle graph one edge removed} we mention that the group associated to the cycle graph on four vertices where we remove one edge also has the $R_{\infty}$--property. Hence, this provides a 7-dimensional example. One can prove (see \cite[Example 8.2.2]{masterthesis}) that this is a sharp bound when considering groups associated to graphs. Moreover, we will show in a forthcoming paper that this is a general lower bound in the sense that there do not exist 2-step nilmanifolds (so also not associated to a graph) in dimensions $\leq 6$ having the $R_{\infty}$--property.
\end{remark}


\section{The second method: simplicial join}
There are several ways to combine graphs. We discuss the \textit{simplicial join} and the \textit{disjoint union}.
\begin{definition}\label{def:disjoint union}
Let $k\in \mathbb{N}_{>1}$ and $\Gamma_i(V_i,E_i)$ (with $i=1,2,\dots,k$) be graphs.
\begin{itemize}
\item The disjoint union $\sqcup_{i=1}^k \Gamma_i$ of the graphs $\Gamma_1,\Gamma_2,\dots, \Gamma_k$ is defined by
\[\bigsqcup_{i=1}^k \Gamma_i\left (\bigsqcup_{i=1}^k V_i, \bigsqcup_{i=1}^k E_i\right ).\]
\item The simplicial join $\ast_{i=1}^k \Gamma_i$ of the graphs $\Gamma_1,\Gamma_2,\dots, \Gamma_k$ is defined by
\[\underset{i=1}{\stackrel{k}{\ast}} \Gamma_i\left (\bigsqcup_{i=1}^k V_i, \bigsqcup_{i=1}^k E_i\cup \{\:x_ix_j\:|\: x_i\in V_i,\: x_j\in V_j \st{and} 1\leq i<j\leq k\}\right ).\]
\end{itemize}
\end{definition}

The group associated to the simplicial join of graphs can be written as a direct product. The next result follows almost immediately by using the definitions.

\begin{lemma}\label{lemma:group associated to simplicial join}
If $\Gamma^{(1)}$ and $\Gamma^{(2)}$ are two undirected simple graphs, then the group $G_{\Gamma^{(1)}\ast \Gamma^{(2)}}$ associated to the simplicial join of $\Gamma^{(1)}$ and $\Gamma^{(2)}$ is isomorphic with the direct product $G_{\Gamma^{(1)}}\times G_{\Gamma^{(2)}}$. This can be generalised to the simplicial join of any finite amount of undirected simple graphs.
\end{lemma}

In order to study the endomorphisms of a direct product, we introduce some notation (which coincides with the notation from \cite{senden_direct}).

\begin{notation}\label{not:induced automorphism components}
Let $G=\bigtimes_{i=1}^k G_i$ be the direct product of $k$ groups $G_i$. For any $i=1,2,\dots,k$ we denote with $\pi_i:G\to G_i$ the canonical projection and with $e_i: G_i\to G$ the canonical inclusion using the direct product. For any endomorphism $\varphi\in \End{G}$ and any $i,j=1,2,\dots,k$, we denote with $\varphi_{ij}:G_j\to G_i$ the morphisms defined by
\[\varphi_{ij}:=\pi_i\circ \varphi\circ e_j:G_j\to G_i.\]
\end{notation}

Note that if $\varphi\in \Aut{G_{\Gamma^{(1)}\ast\Gamma^{(2)}}}$ is an automorphism, then $\varphi_{ii}\in \End{G_{\Gamma^{(i)}}}$ (for $i=1,2$) is not necessarily an automorphism of $G_{\Gamma^{(i)}}$. However, the following lemma tells us something about the images of $\varphi_{ij}$. We refer to \cite[Lemma 2.2]{senden_direct} for a proof.
\begin{lemma}\label{lemma:images induced automorphism components}
Let $\varphi\in \Aut{\bigtimes_{i=1}^k G_i}$ be an automorphism of $\bigtimes_{i=1}^k G_i$, then for all $i=1,2,\dots,k$ it holds that $G_i$ is generated by $\Image{\varphi_{i1}}, \Image{\varphi_{i2}},\dots, \Image{\varphi_{ik}}$.
\end{lemma}

In order to describe the Reidemeister spectrum, we introduce some (obvious) notation.

\begin{notation}\label{not:Reidemeister spectrum product}
Let $M,N\subseteq \mathbb{N}_0\cup \{\infty\}$ be two subsets. We define the product set $M\cdot N$ by
\[M\cdot N:= \{mn\: |\: m\in M,\: n\in N\}.\]
To avoid confusion, we only use this notation when we express the Reidemeister spectrum of a group.
\end{notation}

If $\Gamma(V,E)$ is a graph and $V'\subset V$ a subset of the set of vertices, then we denote with $\Gamma(V')$ the \textit{subgraph induced on $V'$}. This subgraph of $\Gamma$ is defined by means of the vertex set $V'$ and the edge set $\{vw\in E\:|\: v,w\in V'\}$.

A first step in studying the Reidemeister spectrum of the group associated to the simplicial join of graphs is to get rid of the vertices that are connected with all other vertices. Recall that these vertices are precisely contained in the vertex set $V_{n-1}$.
\begin{theorem}\label{thm:graphs abelian parts}
Let $\Gamma(V=\{x_1,x_2,\dots,x_n\},E)$ be an undirected simple graph with $r:=|V_{n-1}|$. If $r<n$, then it holds that
\[\Spec{G_{\Gamma}}=\Spec{\mathbb{Z}^r}\cdot \Spec{G_{\Gamma(V\setminus V_{n-1})}}.\]
\end{theorem}
\begin{proof}
Note that $\Gamma=\Gamma(V_{n-1})\ast \Gamma(V\setminus V_{n-1})$ and thus by Lemma \ref{lemma:group associated to simplicial join}
\[G\cong G_{\Gamma(V_{n-1})}\times G_{\Gamma(V\setminus V_{n-1})}\cong \mathbb{Z}^r\times G_{\Gamma(V\setminus V_{n-1})}.\]
Fix any automorphism $\varphi\in \Aut{\mathbb{Z}^r\times G_{\Gamma(V\setminus V_{n-1})}}$. Since none of the vertices of $\Gamma(V\setminus V_{n-1})$ is connected with all the other vertices of $\Gamma(V\setminus V_{n-1})$, it follows by Lemma \ref{lemma:factor groups Ggamma} that
\[\varphi(\mathbb{Z}^r\times 1)\subseteq \varphi(\mathbb{Z}^r\times Z(G_{\Gamma(V\setminus V_{n-1})}))=Z(\mathbb{Z}^r\times G_{\Gamma(V\setminus V_{n-1})})=\mathbb{Z}^r\times \gamma_2(G_{\Gamma(V\setminus V_{n-1})}).\]
Denote with $(\varphi_{ij})_1$ and $(\varphi_{ij})_2$ (with $i,j=1,2$) the morphisms induced on the first and second factor of the lower central series. Since $\Image{\varphi_{21}}\subseteq \gamma_2(G_{\Gamma(V\setminus V_{n-1})})$ it holds that $\Image{(\varphi_{21})_1}=1$. Hence, $(\varphi_{21})_1$ is the map sending everything to $1\gamma_2(G_{\Gamma(V\setminus V_{n-1})})$. So Lemma \ref{lemma:images induced automorphism components} implies that $(\varphi_{22})_1$ is surjective. Applying Lemma's \ref{lemma:surj ==> aut} and \ref{lemma:phi_1 aut ==> phi aut} yields that $\varphi_{22}$ is an automorphism of $G_{\Gamma(V\setminus V_{n-1})}$.

By taking a particular generating set of $\mathbb{Z}^r$ and $G_{\Gamma(V\setminus V_{n-1})}$ it follows that the matrix of $\varphi_1$ (with respect to this generating set) is of the form
\[\begin{pmatrix}
A & B\\
0 & C
\end{pmatrix}\]
where $A$, $B$ and $C$ are the matrices representing respectively $(\varphi_{11})_1$, $(\varphi_{12})_1$ and $(\varphi_{22})_1$. Since $\varphi_1$ and $(\varphi_{22})_1$ are automorphisms, we can conclude that also $(\varphi_{11})_1$ is an automorphism. Remark that since $\mathbb{Z}^r$ is abelian, it holds that $(\varphi_{11})_1=\varphi_{11}$.

Note that
\[\Image{(\varphi_{12})_2}=\Image{(\varphi_{21})_2}=1\]
and thus $R(\varphi_2)=R((\varphi_{22})_2)$. By Lemma \ref{lemma:Ggamma RN product induced aut} we can conclude that
\[R(\varphi)=R(\varphi_{11}) R((\varphi_{22})_1) R((\varphi_{22})_2)=R(\varphi_{11})R(\varphi_{22}).\]
Since $\varphi_{11}\in\Aut{\mathbb{Z}^r}$ and $\varphi_{22}\in \Aut{G_{\Gamma(V\setminus V_{n-1})}}$ it follows that
\[\Spec{G_{\Gamma}}=\Spec{\mathbb{Z}^r\times G_{\Gamma(V\setminus V_{n-1})}}\subseteq \Spec{\mathbb{Z}^r}\cdot \Spec{G_{\Gamma(V\setminus V_{n-1})}}.\]
The other inclusion is well-known (see e.g. \cite[Corollary 2.6]{senden_direct}).
\end{proof}

\begin{notation}
Using Theorem \ref{thm:graphs abelian parts}, we can restrict ourselves to look at graphs for which none of the vertices is connected (via an edge) with all other vertices. Let $\Gamma$ be such a finite undirected simple graph. Assume that $\Gamma=\ast_{i=1}^k \Gamma^{(i)}$ and that $\Gamma$ cannot be decomposed any further using the simplicial join. We use the superscript ``$(i)$'' to denote similar properties as before, but related to the graph $\Gamma^{(i)}$ (e.g. $V^{(i)}$, $x_1^{(i)}$ and $n^{(i)}$). Applying Lemma \ref{lemma:group associated to simplicial join} yields that $G_{\Gamma}\cong \bigtimes_{i=1}^k G_{\Gamma^{(i)}}$. For any component $i=1,2,\dots,k$ we define the subgroup $H^{(i)}\subseteq G_{\Gamma}$ by
\[H^{(i)}:=\bigtimes_{m=1}^{i-1} \gamma_2(G_{\Gamma^{(m)}})\times G_{\Gamma^{(i)}}\times \bigtimes_{m=i+1}^{k} \gamma_2(G_{\Gamma^{(m)}})=G_{\Gamma^{(i)}}\gamma_2(G_{\Gamma}).\]
\end{notation}

\begin{definition}\label{def:complement graph}
Let $\Gamma(V,E)$ be an undirected simple graph. The complement $\Gamma^{c}$ of the graph $\Gamma$ is defined by
\[\Gamma^{c}(V,\{vw\:|\: v,w\in V, v\neq w \st{and} vw\not\in E\}).\]
\end{definition}

The following lemma will be needed later on.

\begin{lemma}\label{lemma:complement graph}
If $\Gamma=\ast_{i=1}^k \Gamma^{(i)}$ cannot be decomposed any further using the simplicial join, then $\left (\Gamma^{(i)}\right )^c$ is connected (for all $i=1,2,\dots,k$).
\end{lemma}
\begin{proof}
Fix any $i=1,2,\dots,k$ and denote for the sake of simplicity $\Gamma^{(i)}=X(V,E)$. Suppose by contradiction that $X^c$ is not connected. So there exist two subgraphs $X_1(V_1,E_1)$ and $X_2(V_2,E_2)$ of $X^c$ that are not connected with each other (with $V=V_1\sqcup V_2$ and with $E_1\cup E_2$ the edge set of $X^c$). This implies that $V_1,V_2\subset V$ are non-empty sets of vertices such that
\[\{\:x_{i_1}x_{i_2}\:|\: x_{i_1}\in V_1\st{and} x_{i_2}\in V_2\}\subseteq E.\]
Thus, we obtain that
\[X=X(V_1)\ast X(V_2)\]
which contradicts the assumption. Hence, all the complements $\left (\Gamma^{(i)}\right )^c$ are connected.
\end{proof}

Lemma \ref{lemma:complement graph} allows us to describe the automorphisms of $G_{\Gamma}$.

\begin{corollary}\label{cor:simplicial join automorphisms components}
Let $\Gamma(V=\{x_1,x_2,\dots,x_n\},E)$ be an undirected simple graph such that $\Gamma=\ast_{i=1}^k \Gamma^{(i)}$, $|V_{n-1}|=0$ and $\Gamma$ cannot be decomposed any further using the simplicial join. Then, for any automorphism $\varphi\in \Aut{G_{\Gamma}}$ there exists a unique permutation $\sigma\in S_{k}$ such that
\begin{enumerate}
\item[(i)] $\varphi\left (H^{(i)}\right )=H^{(\sigma(i))}$ for any $i=1,2,\dots,k$.
\item[(ii)] The corresponding components are isomorphic, i.e. $\Gamma^{(i)}\cong \Gamma^{(\sigma(i))}$ (for any $i=1,2,\dots,k$).
\end{enumerate}
\end{corollary}
\begin{proof}
Fix any component $i_0\in \{1,2,\dots,k\}$ and denote with $d:=\deg_{\Gamma} (x_1^{(i_0)} )$ the degree (in $\Gamma$) of $x_1^{(i_0)}$. Define for any $i=1,2,\dots,k$ the integer $d^{(i)}$ by
\[d^{(i)}:=d-\sum_{\substack{j=1,\\j\neq i}}^k n^{(j)}=d-n+n^{(i)}.\]
Note that $\deg_{\Gamma^{(i_0)}} (x_1^{(i_0)} )$ equals $d^{(i_0)}$.
\\Since $\varphi$ is an automorphism, it follows by Theorem \ref{thm:degree vertices} that $\varphi (x_1^{(i_0)} )\in H_d\setminus H_{d+1}$. Thus we can fix some component $i_1\in \{1,2,\dots,k\}$ such that $\varphi (x_1^{(i_0)} )$ has a non-zero exponent for some vertex of $V^{(i_1)}$ of degree $d$ in $\Gamma$. By Lemma \ref{lemma:HN centralizer 1 vertex}, it holds that $h (Z_{G_{\Gamma}} (\varphi (x_1^{(i_0)} ) ) )=d+N+1$. However, by applying Lemma \ref{lemma:HN centralizer general element} we obtain (by also using Lemma \ref{lemma:HN properties} $(i)$) that
\begin{align*}
d+N+1&=h (Z_{G_{\Gamma}} (\varphi (x_1^{(i_0)} ) ) )=\sum_{j=1}^k h (Z_{G_{\Gamma^{(j)}}} ( \pi_j( \varphi (x_1^{(i_0)} )) ) )\\
&\leq (d^{(i_1)}+1+N^{(i_1)} )+\sum_{\substack{j=1,\\j\neq i_1}}^k(n^{(j)}+N^{(j)})\\
&=d^{(i_1)}+1-n^{(i_1)}+\sum_{j=1}^k(n^{(j)}+N^{(j)})=d^{(i_1)}+1-n^{(i_1)}+n+N\\
&=d+N+1.
\end{align*}
So equality must hold throughout the calculations and thus for all $j=1,2,\dots,k$ with $j\neq i_1$ this implies that
\[h (Z_{G_{\Gamma^{(j)}}} ( \pi_j( \varphi (x_1^{(i_0)} )) ) )=n^{(j)}+N^{(j)}.\]
Since there are no vertices of degree $n^{(j)}-1$ (in $\Gamma^{(j)}$), Lemma \ref{lemma:HN centralizer general element} implies that $\pi_j( \varphi (x_1^{(i_0)} ))\in\gamma_2(G_{\Gamma^{(j)}})$ for all $j=1,2,\dots,k$ with $j\neq i_1$. So $i_1$ is the unique component such that $\varphi (x_1^{(i_0)} )\in H^{(i_1)}$ and thus we define $\sigma(i_0):=i_1$. Repeating this argument for any $i_0\in \{1,2,\dots,k\}$ yields a unique map $\sigma:\{1,2,\dots,k\}\to\{1,2,\dots,k\}$ such that
\[\varphi(x_1^{(i)})\in H^{(\sigma(i))} \quad \st{for all} i=1,2,\dots,k.\]
Fix any component $i_0$ and any vertex $x_j^{(i_0)}$ of that component (with $j\neq 1$) such that $x_1^{(i_0)}x_j^{(i_0)}\not\in E$. Using the same argument as in the beginning of the proof, we can derive that $\varphi(x_j^{(i_0)})\in H^{(\widetilde{i_0})}$ for some unique component $\widetilde{i_0}$. Since $[x_j^{(i_0)},x_1^{(i_0)}]\neq 1$, it follows that
\[1\neq [\varphi(x_j^{(i_0)}),\varphi(x_1^{(i_0)})]\in [H^{(\widetilde{i_0})}, H^{(\sigma(i_0))}]=\begin{cases}
\gamma_2(G_{\Gamma^{(\sigma(i_0))}}) &\st{if} \widetilde{i_0}= \sigma(i_0)\\
1 &\st{if} \widetilde{i_0}\neq \sigma(i_0)
\end{cases}.\]
Thus we obtain that $\widetilde{i_0}= \sigma(i_0)$ and so $\varphi(x_j^{(i_0)})\in H^{(\sigma(i_0))}$. If we fix any vertex $x_j^{(i_0)}$ (with $j\neq 1$), then by Lemma \ref{lemma:complement graph} there exists a path in $\left (\Gamma^{(i_0)}\right )^c$ that connects $x_j^{(i_0)}$ and $x_1^{(i_0)}$. We can now use the previous argument inductively together with this path to conclude that $\varphi(x_j^{(i_0)})\in H^{(\sigma(i_0))}$.

Hence, we can conclude that $\varphi(H^{(i)})\subseteq H^{(\sigma(i))}$ for all components $i$. By construction of this unique map $\sigma$ and since $\varphi$ is an automorphism, one can derive that $\sigma\in S_k$ and $\varphi(H^{(i)})= H^{(\sigma(i))}$ for all components $i$. The second item follows directly by applying Lemma \ref{lemma:isomorphic graphs}.
\end{proof}

Using this description, we are now able to describe the Reidemeister spectrum of $G_{\Gamma}$. We subdivide this description in two theorems, but we prove them at once.

\begin{theorem}\label{thm:simplicial join isomorphic components}
Let $\Gamma^{(1)}, \Gamma^{(2)},\dots,\Gamma^{(k)}$ be isomorphic finite undirected simple graphs with at least $2$ vertices that cannot be decomposed using the simplicial join. Then
\[\Spec{G_{\ast_{i=1}^k \Gamma^{(i)}}}=\Spec{\bigtimes_{i=1}^k G_{\Gamma^{(i)}}}=\bigcup_{i=1}^{k}\left\{\prod_{m=1}^{i}R_m\: \middle\vert\: R_m\in \Spec{G_{\Gamma^{(1)}}} \right \}.\]
\end{theorem}

If we divide a graph using the simplicial join, then we say that two components \textit{are of the same type} if they are isomorphic. If $s$ is the amount of different types, then we fix some order to be able to address components of type $j$ (with $j=1,2,\dots,s$).

\begin{theorem}\label{thm:simplicial join}
Let $\Gamma(V=\{x_1,x_2,\dots,x_n\},E)$ be a finite undirected simple graph. Assume that $\Gamma(V\setminus V_{n-1})=\ast_{i=1}^k \Gamma^{(i)}$ and $\Gamma(V\setminus V_{n-1})$ cannot be decomposed any further using the simplicial join. Denote with $s$ the number of types of components. Then, the Reidemeister spectrum of $G_{\Gamma}$ is given by
\[\Spec{G_{\Gamma}}=\Spec{\mathbb{Z}^{|V_{n-1}|}}\cdot \prod_{j=1}^s \Spec{\underset{\substack{\text{components }i\\ \text{ of type }j}}{\bigtimes}G_{ \Gamma^{(i)}}}.\]
\end{theorem}
\begin{proof}[Proof of Theorems \ref{thm:simplicial join isomorphic components} and \ref{thm:simplicial join}]
By Theorem \ref{thm:graphs abelian parts} we can assume that $|V_{n-1}|=0$. Fix any automorphism $\varphi\in \Aut{G_{\Gamma}}$. Take $\sigma\in S_k$ as described in Corollary \ref{cor:simplicial join automorphisms components}. Define the map $\overline{\varphi}\in \End{G_{\Gamma}}$ by setting
\[\overline{\varphi}:G_{\Gamma}\to G_{\Gamma}:(g_1,g_2,\dots,g_k)\mapsto \prod_{i=1}^k(1,\dots,1,\underbrace{\varphi_{\sigma(i)i}(g_i)}_{\text{position }\sigma(i)},1,\dots,1).\]
Since the components $\Gamma^{(i)}$ and $\Gamma^{(\sigma(i))}$ are isomorphic, we view $\varphi_{\sigma(i)i}\in \End{G_{\Gamma^{(i)}}}$. Since $\varphi(H^{(i)})=H^{(\sigma(i))}$ (for all $i=1,2,\dots,k$), it follows that if $i_2\neq \sigma(i_1)$ then the induced map $(\varphi_{i_2i_1})_1$ is the zero map. Hence, by using Lemma \ref{lemma:images induced automorphism components}, Lemma \ref{lemma:surj ==> aut} and Lemma \ref{lemma:phi_1 aut ==> phi aut} we obtain that $\varphi_{\sigma(i)i}$ is an automorphism of $G_{\Gamma^{(i)}}$ for all components $i=1,2,\dots,k$ and thus $\overline{\varphi}$ is an automorphism of $G_{\Gamma}$. By using the definition of the automorphism $\overline{\varphi}$, one can derive that $(\overline{\varphi})_1=\varphi_1$ and $(\overline{\varphi})_2=\varphi_2$. Applying Lemma \ref{lemma:Ggamma RN product induced aut} now yields that $R(\varphi)=R(\overline{\varphi})$ and thus it suffices to consider the automorphism $\overline{\varphi}$.

Since $\sigma$ maps components to isomorphic components, we are able to write $\overline{\varphi}=\overline{\varphi}^{(1)}\times \overline{\varphi}^{(2)}\times \dots \times \overline{\varphi}^{(s)}$ for some automorphisms $\overline{\varphi}^{(j)}$ on the groups associated to the graph consisting of the simplicial join of the components of type $j$ (with $j=1,2,\dots,s$). It is known that $R(\overline{\varphi})=\prod_{j=1}^s R(\overline{\varphi}^{(j)})$ (see e.g. \cite[Corollary 2.6]{senden_direct}). So it suffices to prove Theorem \ref{thm:simplicial join isomorphic components}.
\\\\To limit notational complexity, we assume that there are only $k=2$ components. The general case can be proven similarly. We refer the interested reader to \cite[Proposition 8.1.11]{masterthesis}.

Assume that $\Gamma=\Gamma^{(1)}\ast \Gamma^{(2)}$ (and $\Gamma$ cannot be decomposed any further using the simplicial join) where $\Gamma^{(1)}\cong \Gamma^{(2)}$ are isomorphic finite undirected simple graphs with at least $2$ vertices. We argue that $\Spec{G_{\Gamma^{(1)}\ast \Gamma^{(2)}}}=\bigcup_{i=1}^{2}\left\{\prod_{m=1}^{i}R_m\: \middle\vert\: R_m\in \Spec{G_{\Gamma^{(1)}}} \right \}$. First we assume that $\varphi(g_1,1)\in 1\times G_{\Gamma^{(2)}}$ and that $\varphi(1,g_2)\in G_{\Gamma^{(1)}}\times 1$ for all $g_i\in G_{\Gamma^{(i)}}$ (in particular, $\sigma=(1\: 2)\in S_2$). We define a new set of generators for $G_{\Gamma^{(1)}}$ and $G_{\Gamma^{(2)}}$ by setting
\begin{align*}
\tilde{x}_j^{(i)}&:=\begin{cases} x_j^{(1)} &\st{if} i=1\\ \varphi_{21}(x_j^{(1)})&\st{if} i=2\end{cases}\quad \st{for all} j=1,2,\dots,n^{(1)}\\
\tilde{y}_j^{(i)}&:=\begin{cases} y_j^{(1)} &\st{if} i=1\\ \varphi_{21}(y_j^{(1)})&\st{if} i=2\end{cases}\quad \st{for all} j=1,2,\dots,N^{(1)}.
\end{align*}
One can check that this is well-defined since $\varphi$ is an automorphism. Hence, there exists some automorphism $\psi\in \Aut{G_{\Gamma^{(1)}}}$ such that the matrix of $\varphi_i$ (with $i=1,2$) with respect to this new set of generators has the following form
\[\begin{pmatrix}[cc]
0 & A_{i} \\
\mathds{1} & 0
\end{pmatrix}\]
where $A_{i}$ is the matrix of the induced automorphism $\psi_i$. By Lemma \ref{lemma:Ggamma RN product induced aut} we now obtain that
\[R(\varphi_i)=
\left|\det\left(\mathds{1}-\begin{pmatrix}[cc]
0 & A_{i} \\
\mathds{1} & 0
\end{pmatrix}\right)\right|_{\infty}=|\det(\mathds{1}-A_{i})|_{\infty}=R(\psi_i)\]
and hence $R(\varphi)=R(\psi)$.

If on the other hand $\varphi=\varphi_{11}\times\varphi_{22}$ (and thus $\sigma$ is the identity permutation), then $R(\varphi)=R(\varphi_{11})R(\varphi_{22})$ (see e.g. \cite[Corollary 2.6]{senden_direct}). Combining these two cases, we obtain that $\Spec{G_{\Gamma^{(1)}\ast \Gamma^{(2)}}}\subseteq\bigcup_{i=1}^{2}\left\{\prod_{m=1}^{i}R_m\: \middle\vert\: R_m\in \Spec{G_{\Gamma^{(1)}}} \right \}$. One can prove the other inclusion by using the same ideas for constructing the desired automorphisms.
\end{proof}

Using Theorem \ref{thm:simplicial join}, we can now construct more 2-step nilpotent groups having the $R_{\infty}$--property.
\begin{corollary} \label{cor:simplicial join R-infinity}
Let $\Gamma(V,E)$ be a finite undirected simple graph such that $\Gamma=\ast_{i=1}^k \Gamma^{(i)}$ (where $\Gamma$ cannot be decomposed any further). Then the following statements are equivalent
\begin{itemize}
\item[(i)] There exists some $i=1,2,\dots,k$ such that $G_{\Gamma^{(i)}}$ has the $R_{\infty}$--property
\item[(ii)] $G_{\Gamma}$ has the $R_{\infty}$--property
\end{itemize}
\end{corollary}
As a direct consequence of this we now also find the following result.
\begin{corollary}
For any $n\geq 7$ there exists a 2-step nilmanifold $M_n$ of dimension $n$ such that any self-homotopy equivalence of $M_n$ is homotopic to a fixed point free map. 
\end{corollary}
\begin{proof}
Indeed, for $M_7$ we can take the 7-dimensional manifold of Remark~\ref{manifold dim 7}. For any $n>7$, let $k=n-7$ and take 
$M_n=T^k\times M_7$, where $T^k$ is the $k$-dimensional torus. Then, we have that the fundamental group of $M_n$ is 
$\mathbb{Z}^k \times G_\Gamma$, where $\Gamma$ is the cycle graph on four vertices with one edge removed. It follows that the fundamental group of $M_n$ is then the group associated to the simplicial join of $\Gamma$ and $k$ graphs consisting of just one vertex. By the previous corollary, we know that $\mathbb{Z}^k \times G_\Gamma$ has the $R_\infty$--property from which the result follows.
\end{proof}

In order to illustrate Theorem \ref{thm:simplicial join isomorphic components} and Theorem \ref{thm:simplicial join} we give some examples.

\begin{example}\label{ex:simplicial join}
Consider the cycle graph $C_4$ on $4$ vertices in Figure \ref{fig:4V cycle graph}. Note that $C_4$ is the simplicial join of twice the graph with $2$ vertices and no edges. Hence, we can use Theorem \ref{thm:simplicial join isomorphic components} and equation (\ref{eq:Spec N_n,2}) to conclude that
\begin{align*}
\Spec{G_{C_4}}&=\Spec{N_{2,2}\times N_{2,2}}=\bigcup_{i=1}^{2}\left \{\prod_{m=1}^i R_m\: \middle\vert\: R_m\in \Spec{N_{2,2}}\right \}\\
&=\left (2\mathbb{N}_0\cup \{\infty\}\right )\cup\left (4\mathbb{N}_0\cup \{\infty\}\right )
=2\mathbb{N}_0\cup \{\infty\}.
\end{align*}

Let us consider the graph $\Gamma$ in Figure \ref{fig:Graph 4V simplicial join}. Hence, it follows that $V_{n-1}=\{x_1\}$ and thus Theorem \ref{thm:simplicial join} (together with equations (\ref{eq:Spec Z^n}) and (\ref{eq:Spec N_n,2})) yields that
\begin{align*}
\Spec{G_{\Gamma}}&=\Spec{\mathbb{Z}}\cdot \Spec{N_{3,2}}=(\{2,\infty\})\cdot((2\mathbb{N}_0-1)\cup 4\mathbb{N}_0\cup \{\infty\})\\
&=2(2\mathbb{N}_0-1)\cup 8\mathbb{N}_0\cup \{\infty\}.
\end{align*}

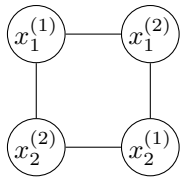
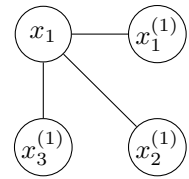
\begin{figure}[H]
     \centering
     \begin{subfigure}[b]{0.4\textwidth}
         \centering
\begin{tikzpicture}
    \tikzstyle{vertex}=
      [circle,draw,minimum size=1.2em,inner sep=0.2];
      
    \node[vertex] (x1) at (0,1.5) {$x_1^{(1)}$};
    \node[vertex] (x2) at (1.5,1.5) {$x_1^{(2)}$};
    \node[vertex] (x3) at (1.5,0) {$x_2^{(1)}$};
    \node[vertex] (x4) at (0,0) {$x_2^{(2)}$};

    \draw (x1)--(x2)--(x3)--(x4)--(x1);
\end{tikzpicture}
         \caption{Cycle graph}
         \label{fig:4V cycle graph}
     \end{subfigure}
     \hfill
     \begin{subfigure}[b]{0.4\textwidth}
         \centering
\begin{tikzpicture}
    \tikzstyle{vertex}=
      [circle,draw,minimum size=2.1em,inner sep=0.2];
      
    \node[vertex] (x1) at (0,1.5) {$x_1$};
    \node[vertex] (x2) at (1.5,1.5) {$x_1^{(1)}$};
    \node[vertex] (x3) at (1.5,0) {$x_2^{(1)}$};
    \node[vertex] (x4) at (0,0) {$x_3^{(1)}$};

    \draw (x4)--(x1)--(x2);
    \draw (x1)--(x3);
\end{tikzpicture}
         \caption{Simplicial join of 1 vertex and 3 vertices}
         \label{fig:Graph 4V simplicial join}
     \end{subfigure}
        \caption{Illustration simplicial join on graphs with $4$ vertices}
        \label{fig:4V simplicial join}
\end{figure}
\end{example}


\section{The third method: connected components}
Let $\Gamma(V,E)$ be a finite undirected simple graph. Denote with $\Gamma^{(0)}$ the induced subgraph $\Gamma(V_0\setminus V_1)$ consisting of vertices of degree zero and with $\Gamma^{(i)}$ (for $i=1,2,\dots,k$) the connected components of $\Gamma(V_1)$. It follows that $\Gamma=\bigsqcup_{i=0}^k \Gamma^{(i)}$. We use the superscript ``$(i)$'' to denote similar properties as before, but related to the graph $\Gamma^{(i)}$. Define for all $i\in \{0,1,\dots,k\}$ the subgroup $H^{(i)}$ of $G_{\Gamma}$ by
\[H^{(i)}:=\left \{\prod_{j=1}^{n^{(i)}} \left (x_j^{(i)}\right )^{z_j^{(i)}}\prod_{l=1}^N y_l^{t_l}\: \middle \vert \:  z_j^{(i)},t_l\in \mathbb{Z}\right \}=G_{\Gamma^{(i)}}\gamma_2(G_\Gamma).\]
As with the simplicial join, it turns out that any automorphism of $G_{\Gamma}$ maps the subgroups $H^{(i)}$ (with $i=1,2,\dots,k$) to such a subgroup associated to an isomorphic component. In order to prove this, we need the following lemma, which studies the centralizers of elements of $G_{\Gamma}$.

\begin{lemma}\label{lemma:HN components graph}
If $\Gamma(V,E)$ is a finite undirected simple graph such that $\Gamma=\bigsqcup_{i=1}^k \Gamma^{(i)}$, then for any $z_j^{(i)},t_l\in \mathbb{Z}$ with some $i_1\neq i_2$ and $j_1,j_2$ such that $z_{j_1}^{(i_1)},z_{j_2}^{(i_2)}\neq 0$, it holds that:
\[Z_{G_{\Gamma}}\left( \prod_{i}\prod_{j} \left (x_j^{(i)}\right )^{z_j^{(i)}}\: \prod_{l=1}^N y_l^{t_l} \right )=\left< \prod_{i}\prod_{j} \left (x_j^{(i)}\right )^{z_j^{(i)}/d}\right>\times \gamma_2(G_{\Gamma})\cong \mathbb{Z}^{N+1}\]
where $d:=\gcd_{i,j} z_j^{(i)}$.
\end{lemma}

\begin{proof}
Without loss of generality we assume that $i_1=1,i_2=2$ and $j_1=j_2=1$ (and thus $z_{1}^{(1)},z_{1}^{(2)}\neq 0$). To make notation more clear, we denote
\[x=\prod_{i}\prod_{j} \left (x_j^{(i)}\right )^{z_j^{(i)}}\: \prod_{l=1}^N y_l^{t_l}.\]
Using the operation in $G_{\Gamma}$, we obtain that
\[Z_{G_{\Gamma}}(x)=\left\{\prod_{i}\prod_{j} \left (x_j^{(i)}\right )^{v_j^{(i)}}\: \prod_{l=1}^N y_l^{s_l}\: \middle\vert \: 
\begin{array}{ll} 
\vspace{5pt}v_j^{(i)},s_l\in \mathbb{Z}& \\
\vspace{5pt} z_{j_1}^{(i_1)}v_{j_2}^{(i_2)}=v_{j_1}^{(i_1)}z_{j_2}^{(i_2)}\: &\st{if} i_1\neq i_2;\: \forall j_1,j_2\\
\vspace{5pt} z_{j_1}^{(i)}v_{j_2}^{(i)}=v_{j_1}^{(i)}z_{j_2}^{(i)}\: &\st{if} x_{j_1}^{(i)}x_{j_2}^{(i)}\not\in E^{(i)}\\
\end{array}\right\}.\]
Since $z_{1}^{(1)},z_{1}^{(2)}\neq 0$, the equations of the form $z_1^{(1)}v_j^{(i)}=v_1^{(1)}z_j^{(i)}$ (with $i=2,3,\dots,k$) and the equations of the form $z_1^{(2)}v_j^{(1)}=v_1^{(2)}z_j^{(1)}$ can only be satisfied if there exists some $\lambda\in \mathbb{Z}$ such that
\[v_j^{(i)}=\lambda \: \frac{z_j^{(i)}}{d}\quad\quad \st{for all} i=1,2,\dots,k\st{and} j=1,2,\dots,n^{(i)}\]
where $d:=\gcd_{i,j} z_{j}^{(i)}$. This solution satisfies all the other conditions in the expression of the centralizer of $x$ and thus we can conclude that
\begin{align*}
Z_{G_{\Gamma}}(x)&=\left\{\left (\prod_{i}\prod_{j} \left (x_j^{(i)}\right )^{z_j^{(i)}/d}\right )^{v_0}\: \prod_{l=1}^N y_l^{s_l}\: \middle \vert \: v_0,s_l\in \mathbb{Z}\right \}=\left< \prod_{i}\prod_{j} \left (x_j^{(i)}\right )^{z_j^{(i)}/d}\right>\times \gamma_2(G_{\Gamma})\cong \mathbb{Z}^{N+1}.
\end{align*}
\end{proof}

Lemma \ref{lemma:HN components graph} allows us to study the automorphisms of $G_{\Gamma}$.

\begin{corollary}\label{cor:aut components graph}
Let $\Gamma(V,E)$ be a finite undirected simple graph. Denote with $\Gamma^{(0)}$ the induced subgraph $\Gamma(V_0\setminus V_1)$ and with $\Gamma^{(i)}$ the connected components of $\Gamma(V_1)$ (where $i=1,2,\dots,k$). Then for any automorphism $\varphi\in \Aut{G_{\Gamma}}$ there exists a unique permutation $\sigma\in S_k$ such that:
\begin{enumerate}
\item[(i)] $\varphi\left (H^{(i)}\right )=H^{(\sigma(i))}$ for any $i=1,2,\dots,k$.
\item[(ii)] The corresponding connected components are isomorphic, i.e. $\Gamma^{(i)}\cong \Gamma^{(\sigma(i))}$ (for any $i=1,2,\dots,k$).
\end{enumerate}
\end{corollary}
\begin{proof}
Fix any component $i_0\in\{1,2,\dots,k\}$. Since $\Gamma^{(i_0)}$ is a subgraph of $\Gamma(V_1)$, it follows by Lemma \ref{lemma:HN centralizer 1 vertex} that
\[h(Z_{G_{\Gamma}}(\varphi(x_1^{(i_0)})))=\deg(x_1^{(i_0)})+N+1>N+1.\]
Lemma \ref{lemma:HN components graph} now implies that $\varphi(x_1^{(i_0)})$ can only have non-zero exponents for vertices from one component. Recall that $H_1$ is a characteristic subgroup (see Lemma \ref{thm:degree vertices}). Hence, we obtain that there is a unique component $\sigma(i_0)\in \{1,2,\dots,k\}$ such that $\varphi(x_1^{(i_0)})\in H^{(\sigma(i_0))}$. Repeating this argument for any $i_0\in \{1,2,\dots,k\}$ yields a unique map $\sigma:\{1,2,\dots,k\}\to\{1,2,\dots,k\}$ such that
\[\varphi(x_1^{(i)})\in H^{(\sigma(i))} \quad \st{for all} i=1,2,\dots,k.\]
Fix any component $i_0\in \{1,2,\dots,k\}$ and some vertex $x_j^{(i_0)}$ of that component such that $x_1^{(i_0)}x_j^{(i_0)}\in E$. Using the same argument as in the beginning of the proof, we can derive that $\varphi(x_j^{(i_0)})\in H^{(\widetilde{i_0})}$ for some unique component $\widetilde{i_0}\in \{1,2,\dots,k\}$. Since $x_1^{(i_0)}$ and $x_j^{(i_0)}$ commute, it holds that also $\varphi(x_1^{(i_0)})$ and $\varphi(x_j^{(i_0)})$ commute. However, since $\varphi(x_1^{(i_0)}),\varphi(x_j^{(i_0)})\not\in \gamma_2(G_{\Gamma})$ and vertices from different components are not connected via an edge (and thus the corresponding group elements do not commute), we obtain that $\widetilde{i_0}=\sigma(i_0)$ and thus $\varphi(x_j^{(i_0)})\in H^{(\sigma(i_0))}$. If $x_j^{(i_0)}$ is any vertex of component $i_0$, then there exists a path in $\Gamma^{(i_0)}$ connecting $x_1^{(i_0)}$ and $x_j^{(i_0)}$. Using the previous argument inductively with this path we can derive that $\varphi(x_j^{(i_0)})\in H^{(\sigma(i_0))}$.

Hence, we can conclude that $\varphi(H^{(i)})\subseteq H^{(\sigma(i))}$ for all components $i$. By construction of this unique map $\sigma$ and since $\varphi$ is an automorphism, one can derive that $\sigma\in S_k$ and $\varphi(H^{(i)})= H^{(\sigma(i))}$ for all components $i\in \{1,2,\dots,k\}$. The second item follows directly by applying Lemma \ref{lemma:isomorphic graphs}.
\end{proof}

In contract to the situation for the simplicial join, it is not possible to give a nice general description of the Reidemeister spectrum of a group $G_\Gamma$ in terms of the Reidemeister spectra of the groups associated to the connected components of $\Gamma$. Nevertheless, in practice Corollary~\ref{cor:aut components graph} is very useful to determine the Reidemeister spectrum in concrete cases.

Indeed, fix any automorphism $\varphi\in \Aut{G_{\Gamma}}$. Corollary \ref{cor:aut components graph} allows us to describe the matrices corresponding to $\varphi_1$ and $\varphi_2$. These matrices will have the following form
\[\begin{pmatrix}[cccc|c]
A^{(1)} & 0 & \hdots & 0 & \ast \\ 
0 & A^{(2)} & \ddots & \vdots & \vdots \\ 
\vdots & \ddots & \ddots & 0 & \vdots \\ 
0 & \hdots & 0 & A^{(s)} & \ast \\ \hline
0 & \hdots & \hdots & 0 & A^{(0)}
\end{pmatrix}
\qt{and}
\begin{pmatrix}[cccc|c|c|c]
(A^{(1)})_2 & 0 & \hdots & 0 & 0 & \ast & \ast \\ 
0 & (A^{(2)})_2 & \ddots & \vdots & \vdots & \vdots & \vdots \\ 
\vdots & \ddots & \ddots & 0 & \vdots & \vdots & \vdots \\ 
0 & \hdots & 0 & (A^{(s)})_2 & 0 & \ast & \ast \\ \hline
0 & \hdots & \hdots & 0 & T & \ast & \ast \\ \hline
0 & \hdots & \hdots & 0 & 0 & T' & \ast \\ \hline
0 & \hdots & \hdots & 0 & 0 & 0 & (A^{(0)})_2
\end{pmatrix}\]
where (as with the simplicial join) we use $s$ to denote the number of types of components. The matrices with a superscript correspond to automorphisms on the induced groups associated to the disjoint union of the components of a particular type.

The matrix $T$ can be described by using blocks corresponding with the commutators between different connected components of $\Gamma(V_1)$. Each column of this block matrix consists of all zeros except at one position. The matrix at this spot can be described by means of the tensor product of two matrices. Similarly, the matrix $T'$ can be described by using blocks corresponding with the commutators between $\Gamma(V_0\setminus V_1)$ and the connected components of $\Gamma(V_1)$. Each column of this block matrix consists of all zeros except at one position which can be described by using the tensor product of some matrix with $A^{(0)}$. For the precise description of the matrices corresponding to $\varphi_1$ and $\varphi_2$, we refer the interested reader to \cite[Application 7.1.6]{masterthesis}. We illustrate this full description by considering two graphs with four vertices.

\begin{example}\label{ex:4V components}
Let us consider the graph $\Gamma$ in Figure \ref{fig:Graph 4V 1E}. Hence, $\Gamma=\Gamma^{(0)}\sqcup \Gamma^{(1)}$ where $\Gamma^{(0)}$ and $\Gamma^{(1)}$ denote the induced subgraphs $\Gamma(V_0)$ and $\Gamma(V_1\setminus V_0)$. Fix any automorphism $\varphi\in \Aut{G_{\Gamma}}$. Using the full description it follows that the matrices of $\varphi_1$ and $\varphi_2$ have the following form
\[\begin{pmatrix}
A^{(1)} & \ast \\ 
0 & A^{(0)}
\end{pmatrix} \qt{and} \begin{pmatrix}
A^{(1)}\otimes A^{(0)} & \ast \\ 
0 & \det(A^{(0)})
\end{pmatrix}\]
where $A^{(1)},A^{(0)}\in \text{GL}_2(\mathbb{Z})$.

Consider the graph $\Gamma$ in Figure \ref{fig:Graph 4V 2E 2C}. Thus $\Gamma=\Gamma^{(1)}\sqcup \Gamma^{(2)}$ where $\Gamma^{(1)}$ and $\Gamma^{(2)}$ denote the two connected components of $\Gamma(V_1\setminus V_0)$. For any automorphism $\varphi\in \Aut{G_{\Gamma}}$ we obtain by Corollary \ref{cor:aut components graph} that either $\varphi(H^{(1)})=H^{(2)}$ and $\varphi(H^{(2)})=H^{(1)}$ or that $\varphi(H^{(i)})=H^{(i)}$ (for $i=1,2$). In the first case, the matrices of $\varphi_1$ and $\varphi_2$ have the following form
\[\begin{pmatrix}
0 & A^{(1)}_2 \\
A^{(1)}_1 & 0
\end{pmatrix} \qt{and} \begin{pmatrix}
-(A^{(1)}_2\otimes A^{(1)}_1)
\end{pmatrix}\]
where $A^{(1)}_1,A^{(1)}_2\in \text{GL}_2(\mathbb{Z})$. In the other case, the matrices of $\varphi_1$ and $\varphi_2$ have the following form
\[\begin{pmatrix}
A^{(1)}_1 & 0 \\ 
0 & A^{(1)}_2
\end{pmatrix} \qt{and} \begin{pmatrix}
A^{(1)}_1\otimes A^{(1)}_2
\end{pmatrix}\]
where $A^{(1)}_1,A^{(1)}_2\in \text{GL}_2(\mathbb{Z})$.

For the details of these two examples we refer the reader to \cite[Theorem 7.2.2 and 7.2.3]{masterthesis}.
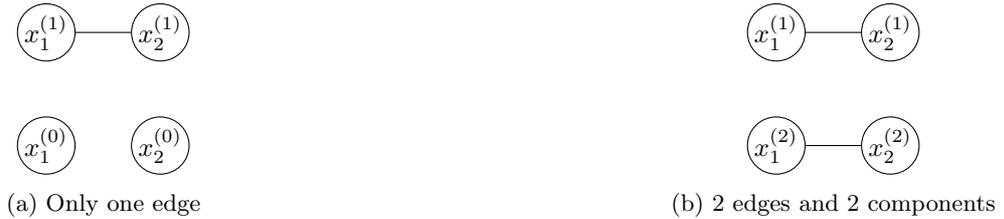
\begin{figure}[H]
     \centering
     \begin{subfigure}[b]{0.4\textwidth}
         \centering
\begin{tikzpicture}
    \tikzstyle{vertex}=
      [circle,draw,minimum size=1.2em,inner sep=0.2];
      
    \node[vertex] (x1) at (0,1.5) {$x_1^{(1)}$};
    \node[vertex] (x2) at (1.5,1.5) {$x_2^{(1)}$};
    \node[vertex] (x3) at (1.5,0) {$x_2^{(0)}$};
    \node[vertex] (x4) at (0,0) {$x_1^{(0)}$};

    \draw (x1)--(x2);
\end{tikzpicture}
         \caption{Only one edge}
         \label{fig:Graph 4V 1E}
     \end{subfigure}
     \hfill
     \begin{subfigure}[b]{0.4\textwidth}
         \centering
\begin{tikzpicture}
    \tikzstyle{vertex}=
      [circle,draw,minimum size=1.2em,inner sep=0.2];
      
    \node[vertex] (x1) at (0,1.5) {$x_1^{(1)}$};
    \node[vertex] (x2) at (1.5,1.5) {$x_2^{(1)}$};
    \node[vertex] (x3) at (1.5,0) {$x_2^{(2)}$};
    \node[vertex] (x4) at (0,0) {$x_1^{(2)}$};

    \draw (x1)--(x2);
    \draw (x3)--(x4);
\end{tikzpicture}
         \caption{$2$ edges and $2$ components}
         \label{fig:Graph 4V 2E 2C}
     \end{subfigure}
        \caption{Illustration disjoint union on graphs with $4$ vertices}
        \label{fig:4V components}
\end{figure}
\end{example}


\section{Examples}

We illustrate the results from the previous sections by determining the Reidemeister spectrum of some families of groups associated to graphs. We first introduce some notation.
\begin{notation}
Let $\Gamma(\{x_1,x_2,\dots,x_n\},E)$ be a finite undirected simple graph and denote with $\boldsymbol{\lambda}=\begin{pmatrix}
\lambda_1 & \lambda_2 & \hdots & \lambda_n
\end{pmatrix}^T\in \mathbb{Z}^{n\times 1}$ some vector of integers. We introduce the notation $(x_1,\dots,x_n)^{\boldsymbol{\lambda}}$ to denote $x_1^{\lambda_1}x_2^{\lambda_2}\dots x_n^{\lambda_n}$. If $A\in \mathbb{Z}^{n\times n}$, then we denote by $A_{:i}$ (respectively $A_{i:}$) the $i$-th column (respectively row) of $A$ (with $i=1,2,\dots,n$).
\end{notation}
 
\subsection{Disjoint union of a complete graph and an isolated vertex}

Denote with $\Gamma_n$ (for $n\in \mathbb{N}_{>2}$) the disjoint union of the complete graph on $n-1$ vertices and an isolated vertex. We assume that $x_n$ is the isolated vertex (and thus $\Gamma(\{x_1,x_2,\dots,x_{n-1}\})$ is a complete graph). Denote with $y_i:=[x_n,x_i]$ (for $i=1,2,\dots,n-1$) the commutators of $G_{\Gamma_n}$.

\begin{theorem}\label{thm:RS isolated + complete}
The Reidemeister spectrum of the groups associated to $\Gamma_n$ (for $n\in \mathbb{N}_{>2}$) is given by
\[\Spec{G_{\Gamma_n}}=\begin{cases}
2\mathbb{N}_0^2\: \cup\: 2|\mathbb{N}^2-4|_{\infty}\cup \{\infty\} & \text{if } n=3\\
2(2\mathbb{N}_0-1)\cup 8\mathbb{N}_0\cup \{\infty\} & \text{if } n\geq 4
\end{cases}\]
where $\mathbb{N}^2$ (respectively $\mathbb{N}_0^2$) denotes the squares (respectively non-zero squares) of integers.
\end{theorem}

\begin{proof}
We use a similar approach as in \cite[section 4]{dekimpe2020} where they determine the Reidemeister spectrum of $N_{r,2}$ for $r\in \mathbb{N}_{\geq 2}$.

Fix any automorphism $\varphi\in \Aut{G_{\Gamma_n}}$. Corollary \ref{cor:aut components graph} implies that the matrix of $\varphi_1$ has the following form
\[\begin{pmatrix}[cc]
A & *\\
  0 & \alpha
\end{pmatrix}\]
where $\ast\in \mathbb{Z}^{(n-1)\times 1}$, $A\in \text{GL}_{n-1}(\mathbb{Z})$ and $\alpha\in \{-1,1\}$. Note that by Lemma \ref{lemma:bilinearity [.,.] 2-nilpotent} for any $i=1,2,\dots, n-1$ it holds that
\[\varphi(y_i)=[\varphi(x_n),\varphi(x_i)]=[x_n^{\alpha},(x_1,\dots,x_{n-1})^{A_{:i}}]=(y_1,\dots,y_{n-1})^{\alpha A_{:i}}\]
and thus the matrix of $\varphi_2$ is equal to $\alpha A$. Using Lemma \ref{lemma:Ggamma RN product induced aut} it follows that
\[R(\varphi)=|\det(\mathds{1}_{n-1}-A)|_{\infty}\: |1-\alpha|_{\infty}\: |\det(\mathds{1}_{n-1}-\alpha A)|_{\infty}.\]
If we assume that $R(\varphi)<\infty$, then we obtain that $\alpha=-1$. We denote with $p_A\in \mathbb{Z}[x]$ the characteristic polynomial of the matrix $A$. Hence, we get that $R(\varphi)=2|p_A(1)p_A(-1)|_{\infty}$.

To any matrix $A\in \text{GL}_{n-1}(\mathbb{Z})$ we can associate an automorphism of $G_{\Gamma_n}$ with Reidemeister number $2|p_A(1)p_A(-1)|_{\infty}$. Indeed, fix any matrix $A\in \text{GL}_{n-1}(\mathbb{Z})$. Define the map $\varphi:G_{\Gamma_n}\to G_{\Gamma_n}$ by setting (with $i=1,2,\dots,n-1$):
\[\begin{cases}
\varphi(x_i)=(x_1,\dots,x_{n-1})^{A_{:i}}\\
\varphi(x_n)=x_n^{-1}\\
\varphi(y_i)=(y_1,\dots,y_{n-1})^{-A_{:i}}
\end{cases}\]
and extending it to $G_{\Gamma_n}$. One can check that $\varphi\in \Aut{G_{\Gamma_n}}$ and that $R(\varphi)=2|p_A(1)p_A(-1)|_{\infty}$. For any monic polynomial $p(x)=x^{n-1}+a_{n-2}x^{n-2}+\dots+a_0\in \mathbb{Z}[x]$ (with $a_0=\pm 1$) we consider the \textit{companion matrix} $C_p$ of the polynomial $p$, i.e. the matrix defined by 
\[C_p=\begin{pmatrix}[ccc|c]
0 & \hdots & 0 & -a_0 \\ \hline
  &   &   & -a_1 \\ 
  & \mathds{1}_{n-2} &   & \vdots \\ 
  &   &   & -a_{n-2}
\end{pmatrix}.\]
Note that $C_p\in \text{GL}_{n-1}(\mathbb{Z})$ and its characteristic polynomial is equal to $p$. Thus using the matrix $C_p$ and the above argument, it follows that there exists an automorphism of $G_{\Gamma_n}$ with Reidemeister number $2|p(1)p(-1)|_{\infty}$. Hence, we obtain that
\[\Spec{G_{\Gamma_n}}=\{2|p(1)p(-1)|_{\infty}\: |\: p(x)\in \mathbb{Z}[x] \st{is monic,} \deg(p)=n-1 \st{and} p(0)=\pm 1\}.\]
Fix any $p(x)=x^{n-1}+a_{n-2}x^{n-2}+\dots+a_0\in \mathbb{Z}[x]$ (with $a_0=\pm 1$). If $n=3$, then
\[|p(1)p(-1)|_{\infty}=2|(1+a_1+a_0)(1-a_1+a_0)|_{\infty}=\begin{cases}
2|a_1^2|_{\infty} &\st{if} a_0=-1\\
2|4-a_1^2|_{\infty} &\st{if} a_0=1
\end{cases}.\]
Using this, it indeed follows that
\[\Spec{G_{\Gamma_3}}=2\mathbb{N}_0^2\: \cup\: 2|\mathbb{N}^2-4|_{\infty}\: \cup \: \{\infty\}.\]
If $n\geq 4$, then one can check that
\[2|p(1)p(-1)|_{\infty}=\begin{cases}
2\left|\left (1+\sum_{i=0}^{m-1}a_{2i}\right )^2-\; \left (\sum_{i=1}^m a_{2i-1}\right )^2\;\; \right|_{\infty} &\st{if} n=2m+1 \vspace{5pt}\\
2\left|\left (\sum_{i=0}^{m-1}a_{2i}\right )^2-\left (1+\sum_{i=1}^{m-1} a_{2i-1}\right )^2\right|_{\infty} &\st{if} n=2m
\end{cases}.\]
So in both cases the Reidemeister spectrum is two times the difference of two squares. Note that the difference of two squares is always a multiple of four or an odd number. Hence, for any $n\in \mathbb{N}_{\geq 4}$ it holds that
\[\Spec{G_{\Gamma_n}}\subseteq 2(2\mathbb{N}_0-1)\cup 8\mathbb{N}_0\cup \{\infty\}.\]
To prove equality, it suffices to define (for any $n\in \mathbb{N}_{\geq 4}$ and $k\in \mathbb{N}_0$) the polynomials
\begin{align*}
q_k(x)&=\begin{cases}
x^{2m}+(k-2)x^2+(k-1)x+1 &\st{if} n=2m+1\\
x^{2m-1}+(k-1)x^2+(k-2)x+1 &\st{if} n=2m
\end{cases}\\
r_k(x)&=\begin{cases}
x^{2m}+(k-1)x^2+(k-1)x+1 &\st{if} n=2m+1\\
x^{2m-1}+kx^2+(k-2)x+1 &\st{if} n=2m
\end{cases}
\end{align*}
and note that
\[2|q_k(1)q_k(-1)|_{\infty}=2(2k-1) \qt{and} 2|r_k(1)r_k(-1)|_{\infty}=8k.\]
\end{proof}

\subsection{Graphs with four vertices and one or two disjoint edges}
In this section we again consider the two graphs from Example \ref{ex:4V components} (see Figure \ref{fig:4V components}). We include these two graphs to illustrate that not all Reidemeister spectra can be described using short and easy expressions. In Example \ref{ex:4V components} we described the matrices of $\varphi_1$ and $\varphi_2$ (for any $\varphi\in \Aut{G_{\Gamma}}$). Using these matrices, Lemma \ref{lemma:Ggamma RN product induced aut} gives us an expression of the Reidemeister number $R(\varphi)$. However, this expression depends on the tensor product of two invertible $2\times 2$ matrices over $\mathbb{Z}$. The next lemma follows by some easy calculations. We refer to \cite[Lemma 7.2.1]{masterthesis} for a detailed proof.

\begin{lemma}\label{lemma:calculations tensor product}
Let $A,B\in \text{GL}_{2}(\mathbb{Z})$ be invertible matrices and $\epsilon=\pm 1$ and let $t_A=\text{Tr}(A)$, $t_B=\text{Tr}(B)$, $d_A=\text{det}(A)$ and $d_B=\text{det}(B)$. Then
\begin{align*}
\det(\mathds{1}_2-A)=1-t_A+d_A&=\begin{cases}2-t_A &\st{if} d_A=1\\-t_A &\st{if} d_A=-1\end{cases}\\
\det(\mathds{1}_4-\epsilon \: A\otimes B)&=\begin{cases}
(t_B-\epsilon \: t_A)^2 &\st{if} d_A=d_B=1\\
-(t_B+\epsilon \: t_A)^2 &\st{if} d_A=d_B=-1\\
-(t_B^2-t_A^2-4) &\st{if} d_A=-1\st{and} d_B=1\\
t_B^2-t_A^2+4 &\st{if} d_A=1\st{and} d_B=-1.
\end{cases}
\end{align*}
\end{lemma}

The expressions from Lemma \ref{lemma:calculations tensor product} together with the description in Example \ref{ex:4V components} allow us to express the Reidemeister number $R(\varphi)$ in terms of two integers. To prove that any such expression is contained in the spectrum, one can use matrices of the form $\begin{pmatrix}
0 & 1 \\ 
1 & m
\end{pmatrix}\in \text{GL}_{2}(\mathbb{Z})$ (with $m\in \mathbb{Z}$) to construct automorphisms having these specified Reidemeister numbers. We omit the details and refer the interested reader to \cite[Theorems 7.2.2 and 7.2.3]{masterthesis}.

\begin{theorem}\label{thm:4V 1E or 2E}
If $\Gamma$ is the undirected simple graph with four vertices and one edge (see Figure \ref{fig:Graph 4V 1E}), then 
\[\Spec{G_{\Gamma}}= \{2|nm(n+m)^2|_{\infty},\: 2|nm(n^2-m^2-4m)|_{\infty} \st{with}m,n\in \mathbb{Z}\}.\]
If $\Gamma$ is the undirected simple graph with four vertices and two disjoint edges (see Figure \ref{fig:Graph 4V 2E 2C}), then 
\[\Spec{G_{\Gamma}}=
\mathbb{N}_0^3\cup \left \{
|nm(n+m)^2|_{\infty},\: |nm(n^2-m^2-4m)|_{\infty},\: |(n-2)(n+2)^2|_{\infty}
\st{with}m,n\in \mathbb{Z}\right \}.\]
\end{theorem}

\subsection{Cycle graphs}
In this section we consider the cycle graphs $C_n$. It turns out that for $n\in \mathbb{N}_{\geq 5}$ the associated group $G_{C_n}$ has the $R_{\infty}$--property. We present the ideas behind this claim, but for the details and the proofs in this section, we refer to \cite[Section 9.1]{masterthesis}.

For any $n\in \mathbb{N}_{\geq 3}$ we define the cycle graph $C_n$ by
\[C_n(\{x_1,x_2,\dots,x_n\},\{x_1x_2,x_2x_3,\dots,x_{n-1}x_n,x_nx_1\}).\] 
In order to describe the Reidemeister spectrum of $G_{C_n}$ we need to understand the automorphisms of $G_{C_n}$. As before, one can study the Hirsch number of the centralizers of elements of $G_{C_n}$ and derive the next lemma.

\begin{lemma}\label{lemma:cycle graph HN}
Let $n\geq 5$ be an integer. For any $z_i,t_l\in \mathbb{Z}$ (for $i=1,2,\dots,n$ and $l=1,2,\dots, N$) it holds that 
\[h\left(Z_{G_{C_n}}\left( \prod_{i=1}^n x_i^{z_i}\: \prod_{l=1}^N y_l^{t_l}\right) \right )<N+3\]
if there are two different indices $i_1,i_2\in \{1,2,\dots,n\}$ such that $z_{i_1},z_{i_2}\neq 0$.
\end{lemma}

By Lemma \ref{lemma:HN centralizer 1 vertex} it holds that $h(Z_{G_{C_n}}(\varphi(x_{i})))=N+3$ for any $i=1,2,\dots,n$ and any $\varphi\in \Aut{G_{C_n}}$. Hence, Lemma \ref{lemma:cycle graph HN} yields that $\varphi_1(x_i\gamma_2(G_{C_n}))\in \langle x_{\sigma(i)}\gamma_2(G_{C_n})\rangle$ for some $\sigma(i)=1,2,\dots,n$. Based on this argument, one can proof the next corollary.

\begin{corollary}\label{cor:cycle graph aut}
Let $n\geq 5$ be an integer and $\varphi\in \Aut{G_{C_n}}$ an automorphism of $G_{C_n}$. Then there exist a permutation $\sigma\in S_n$ and $\epsilon_i=\pm 1$ (for all $i=1,2,\dots,n$) such that
\[\varphi_1(x_i\gamma_2(G_{C_n}))=x_{\sigma(i)}^{\epsilon_i}\gamma_2(G_{C_n})\quad \st{for all} i=1,2,\dots,n.\]
Moreover, the permutation $\sigma$ belongs to the dihedral group $D_n$ of order $2n$. In particular, $\sigma$ consists of a rotation or a reflection of the graph $C_n$.
\end{corollary}

Corollary \ref{cor:cycle graph aut} gives us enough information on the automorphisms to prove that $G_{C_n}$ has the $R_{\infty}$--property for $n\geq 5$. This proof uses Proposition \ref{prop:equivalent statements R-inf} by describing an eigenvector with eigenvalue $1$ for $\varphi_1$ or $\varphi_2$. Note that $G_{C_3}\cong \mathbb{Z}^3$ since $C_3$ is the complete graph on $3$ vertices. We already discussed the cycle graph $C_4$ in Example \ref{ex:simplicial join}. This leads to the next result.

\begin{theorem}\label{thm:cycle graph}
The Reidemeister spectrum of the group $G_{C_n}$ associated to the cycle graph $C_n$ on $n$ vertices is given by
\[\Spec{G_{C_n}}=\begin{cases}
\mathbb{N}_0\cup \{\infty\} &\st{if} n=3\\
2\mathbb{N}_0\cup \{\infty\} &\st{if} n=4\\
\{\infty\} &\st{if} n\geq 5
\end{cases}.\]
\end{theorem}

\begin{remark}\label{rem:cycle graph one edge removed}
Using similar ideas, one can consider the cycle graph $P_n$ on $n$ vertices with one edge removed (also known as the path graph on $n$ vertices). For the details we refer to \cite[Section 9.2]{masterthesis}. The Reidemeister spectrum of $G_{P_n}$ is given by
\[\Spec{G_{P_n}}=\begin{cases}
4\mathbb{N}_0\cup \{\infty\} &\st{if} n=3\\
\{\infty\} &\st{if} n\geq 4
\end{cases}.\]
In \cite[Example 4.1]{gonccalves2009twisted} D.~Gonçalves and P.~Wong prove that $G_{P_4}$ has the $R_{\infty}$--property. We linked the group to the path graph on four vertices to better understand its structure. The methods in this paper provide a more general framework to consider similar groups and allow for a more elegant proof of this result. Moreover, we presented new examples of finitely generated torsion-free $2$-step nilpotent groups that are associated to a graph and have the $R_{\infty}$--property.
\end{remark}

\subsection{Reidemeister spectrum of groups associated to graphs with at most 4 vertices}
Looking back at the different methods we developed, we are able to describe the Reidemeister spectrum of the groups associated to the graphs with at most four vertices. The result is summarised in Table \ref{tab:Reidemeister 3 vertices} and Table \ref{tab:Reidemeister 4 vertices}.

\begin{table}[H]
\centering
\begin{tabular}{ >{\centering\arraybackslash} m{1.4cm}  >{\centering\arraybackslash} m{5.9cm} || >{\centering\arraybackslash} m{1.4cm}  >{\centering\arraybackslash} m{5.9cm} }
Graph $\Gamma$ & Reidemeister spectrum $\Spec{G_{\Gamma}}$ & Graph $\Gamma$ & Reidemeister spectrum $\Spec{G_{\Gamma}}$\\ 
\hline
\vspace{10pt}
\begin{tikzpicture}
    \tikzstyle{vertex}=
      [circle,draw,minimum size=1.2em,inner sep=0.2];
      
    \node[vertex] (x1) at (0,0) {$x_1$};

    \draw (x1);
\end{tikzpicture} \vspace{10pt}& $\{2,\infty\}$
& & \\
\hline
\hline
\vspace{10pt}
\begin{tikzpicture}
    \tikzstyle{vertex}=
      [circle,draw,minimum size=1.2em,inner sep=0.2];
      
    \node[vertex] (x1) at (0,1.3) {$x_1$};
    \node[vertex] (x2) at (1.3,0) {$x_2$};

\end{tikzpicture}\vspace{10pt}& $2\mathbb{N}_0\: \cup\: \{\infty\}$
&
\vspace{10pt}
\begin{tikzpicture}
    \tikzstyle{vertex}=
      [circle,draw,minimum size=1.2em,inner sep=0.2];
      
    \node[vertex] (x1) at (0,1.3) {$x_1$};
    \node[vertex] (x2) at (1.3,0) {$x_2$};

    \draw (x1)--(x2);
\end{tikzpicture}\vspace{10pt} & $\mathbb{N}_0\: \cup\: \{\infty\}$\\ 
\hline
\hline
\vspace{10pt}
\begin{tikzpicture}
    \tikzstyle{vertex}=
      [circle,draw,minimum size=1.2em,inner sep=0.2];
      
    \node[vertex] (x1) at (0,1.3) {$x_1$};
    \node[vertex] (x2) at (1.3,1.3) {$x_2$};
    \node[vertex] (x3) at (1.3,0) {$x_3$};

\end{tikzpicture} \vspace{10pt}& $(2\mathbb{N}_0-1)\: \cup\: 4\mathbb{N}_0\: \cup\: \{\infty\}$ 
&
\vspace{10pt}
\begin{tikzpicture}
    \tikzstyle{vertex}=
      [circle,draw,minimum size=1.2em,inner sep=0.2];
      
    \node[vertex] (x1) at (0,1.3) {$x_1$};
    \node[vertex] (x2) at (1.3,1.3) {$x_2$};
    \node[vertex] (x3) at (1.3,0) {$x_3$};

    \draw (x1)--(x2);
\end{tikzpicture}\vspace{10pt} & $2\mathbb{N}_0^2\: \cup\: 2|\mathbb{N}^2-4|_{\infty}\: \cup \: \{\infty\}$\\
\hline
\vspace{10pt}
\begin{tikzpicture}
    \tikzstyle{vertex}=
      [circle,draw,minimum size=1.2em,inner sep=0.2];
      
    \node[vertex] (x1) at (0,1.3) {$x_1$};
    \node[vertex] (x2) at (1.3,1.3) {$x_2$};
    \node[vertex] (x3) at (1.3,0) {$x_3$};

    \draw (x1)--(x2)--(x3);
\end{tikzpicture} \vspace{10pt}& $4\mathbb{N}_0\: \cup\: \{\infty\}$
& 
\vspace{10pt}
\begin{tikzpicture}
    \tikzstyle{vertex}=
      [circle,draw,minimum size=1.2em,inner sep=0.2];
      
    \node[vertex] (x1) at (0,1.3) {$x_1$};
    \node[vertex] (x2) at (1.3,1.3) {$x_2$};
    \node[vertex] (x3) at (1.3,0) {$x_3$};

    \draw (x1)--(x2)--(x3)--(x1);
\end{tikzpicture} \vspace{10pt}& $\mathbb{N}_0\: \cup\: \{\infty\}$ \\ 
\end{tabular}
\caption{Reidemeister spectrum of groups associated to graphs with at most 3 vertices.}
\label{tab:Reidemeister 3 vertices}
\end{table}

\begin{table}[H]
\centering
\begin{tabular}{ >{\centering\arraybackslash} m{1.4cm}  >{\centering\arraybackslash} m{5.9cm} || >{\centering\arraybackslash} m{1.4cm}  >{\centering\arraybackslash} m{5.9cm} }
Graph $\Gamma$ & Reidemeister spectrum $\Spec{G_{\Gamma}}$ & Graph $\Gamma$ & Reidemeister spectrum $\Spec{G_{\Gamma}}$\\ 
\hline
\vspace{10pt}
\begin{tikzpicture}
    \tikzstyle{vertex}=
      [circle,draw,minimum size=1.2em,inner sep=0.2];
      
    \node[vertex] (x1) at (0,1.3) {$x_1$};
    \node[vertex] (x2) at (1.3,1.3) {$x_2$};
    \node[vertex] (x3) at (1.3,0) {$x_3$};
    \node[vertex] (x4) at (0,0) {$x_4$};
\end{tikzpicture}\vspace{10pt} & $\mathbb{N}_0\: \cup\: \{\infty\}$
&
\vspace{10pt}
\begin{tikzpicture}
    \tikzstyle{vertex}=
      [circle,draw,minimum size=1.2em,inner sep=0.2];
      
    \node[vertex] (x1) at (0,1.3) {$x_1$};
    \node[vertex] (x2) at (1.3,1.3) {$x_2$};
    \node[vertex] (x3) at (1.3,0) {$x_3$};
    \node[vertex] (x4) at (0,0) {$x_4$};

    \draw (x1)--(x2);
\end{tikzpicture}\vspace{10pt} & $\left \{\begin{array}{c}
2|nm(n+m)^2|_{\infty},\\ 2|nm(n^2-m^2-4m)|_{\infty} \\\text{with }m,n\in \mathbb{Z}
\end{array}\right \}$\vspace{3pt}\\
\hline
\vspace{10pt}
\begin{tikzpicture}
    \tikzstyle{vertex}=
      [circle,draw,minimum size=1.2em,inner sep=0.2];
      
    \node[vertex] (x1) at (0,1.3) {$x_1$};
    \node[vertex] (x2) at (1.3,1.3) {$x_2$};
    \node[vertex] (x3) at (1.3,0) {$x_3$};
    \node[vertex] (x4) at (0,0) {$x_4$};

    \draw (x1)--(x2);
    \draw (x3)--(x4);
\end{tikzpicture}\vspace{10pt} & \makecell{$\mathbb{N}_0^3\: \cup$\\$\left \{\begin{array}{c}
|nm(n+m)^2|_{\infty},\\ |nm(n^2-m^2-4m)|_{\infty},\\ |(n-2)(n+2)^2|_{\infty}\\\text{with }m,n\in \mathbb{Z}
\end{array}\right \}$}\vspace{3pt}
&
\vspace{10pt}
\begin{tikzpicture}
    \tikzstyle{vertex}=
      [circle,draw,minimum size=1.2em,inner sep=0.2];
      
    \node[vertex] (x1) at (0,1.3) {$x_1$};
    \node[vertex] (x2) at (1.3,1.3) {$x_2$};
    \node[vertex] (x3) at (1.3,0) {$x_3$};
    \node[vertex] (x4) at (0,0) {$x_4$};

    \draw (x1)--(x2)--(x3);
\end{tikzpicture}\vspace{10pt} & $\{\infty\}$\\ 
\hline
\vspace{10pt}
\begin{tikzpicture}
    \tikzstyle{vertex}=
      [circle,draw,minimum size=1.2em,inner sep=0.2];
      
    \node[vertex] (x1) at (0,1.3) {$x_1$};
    \node[vertex] (x2) at (1.3,1.3) {$x_2$};
    \node[vertex] (x3) at (1.3,0) {$x_3$};
    \node[vertex] (x4) at (0,0) {$x_4$};

    \draw (x1)--(x2)--(x3)--(x1);
\end{tikzpicture}\vspace{10pt} & $2(2\mathbb{N}_0-1)\: \cup\: 8\mathbb{N}_0\: \cup\: \{\infty\}$
&
\vspace{10pt}
\begin{tikzpicture}
    \tikzstyle{vertex}=
      [circle,draw,minimum size=1.2em,inner sep=0.2];
      
    \node[vertex] (x1) at (0,1.3) {$x_1$};
    \node[vertex] (x2) at (1.3,1.3) {$x_2$};
    \node[vertex] (x3) at (1.3,0) {$x_3$};
    \node[vertex] (x4) at (0,0) {$x_4$};

    \draw (x4)--(x1)--(x2);
    \draw (x1)--(x3);
\end{tikzpicture}\vspace{10pt} & $2(2\mathbb{N}_0-1)\:\cup\: 8\mathbb{N}_0\: \cup\: \{\infty\}$\\
\hline
\vspace{10pt}
\begin{tikzpicture}
    \tikzstyle{vertex}=
      [circle,draw,minimum size=1.2em,inner sep=0.2];
      
    \node[vertex] (x1) at (0,1.3) {$x_1$};
    \node[vertex] (x2) at (1.3,1.3) {$x_2$};
    \node[vertex] (x3) at (1.3,0) {$x_3$};
    \node[vertex] (x4) at (0,0) {$x_4$};

    \draw (x4)--(x1)--(x2)--(x3);
\end{tikzpicture}\vspace{10pt} & $\{\infty\}$
&
\vspace{10pt}
\begin{tikzpicture}
    \tikzstyle{vertex}=
      [circle,draw,minimum size=1.2em,inner sep=0.2];
      
    \node[vertex] (x1) at (0,1.3) {$x_1$};
    \node[vertex] (x2) at (1.3,1.3) {$x_2$};
    \node[vertex] (x3) at (1.3,0) {$x_3$};
    \node[vertex] (x4) at (0,0) {$x_4$};

    \draw (x4)--(x1)--(x2);
    \draw (x1)--(x3)--(x2);
\end{tikzpicture}\vspace{10pt} & $4\mathbb{N}_0^2\: \cup\: 4|\mathbb{N}^2-4|_{\infty}\: \cup \: \{\infty\}$\\
\hline
\vspace{10pt}
\begin{tikzpicture}
    \tikzstyle{vertex}=
      [circle,draw,minimum size=1.2em,inner sep=0.2];
      
    \node[vertex] (x1) at (0,1.3) {$x_1$};
    \node[vertex] (x2) at (1.3,1.3) {$x_2$};
    \node[vertex] (x3) at (1.3,0) {$x_3$};
    \node[vertex] (x4) at (0,0) {$x_4$};

    \draw (x1)--(x2)--(x3)--(x4)--(x1);
\end{tikzpicture}\vspace{10pt} & $2\mathbb{N}_0\: \cup\: \{\infty\}$
&
\vspace{10pt}
\begin{tikzpicture}
    \tikzstyle{vertex}=
      [circle,draw,minimum size=1.2em,inner sep=0.2];
      
    \node[vertex] (x1) at (0,1.3) {$x_1$};
    \node[vertex] (x2) at (1.3,1.3) {$x_2$};
    \node[vertex] (x3) at (1.3,0) {$x_3$};
    \node[vertex] (x4) at (0,0) {$x_4$};

    \draw (x4)--(x1)--(x2)--(x3);
    \draw (x1)--(x3);
    \draw (x2)--(x4);
\end{tikzpicture}\vspace{10pt} & $2\mathbb{N}_0\: \cup\:\{\infty\}$\\
\hline
\vspace{10pt}
\begin{tikzpicture}
    \tikzstyle{vertex}=
      [circle,draw,minimum size=1.2em,inner sep=0.2];
      
    \node[vertex] (x1) at (0,1.3) {$x_1$};
    \node[vertex] (x2) at (1.3,1.3) {$x_2$};
    \node[vertex] (x3) at (1.3,0) {$x_3$};
    \node[vertex] (x4) at (0,0) {$x_4$};

    \draw (x1)--(x2)--(x3)--(x4)--(x1)--(x3);
    \draw (x2)--(x4);
\end{tikzpicture}\vspace{10pt} & $\mathbb{N}_0\: \cup\: \{\infty\}$
& \\

\end{tabular}
\caption{Reidemeister spectrum of groups associated to graphs with $4$ vertices.}
\label{tab:Reidemeister 4 vertices}
\end{table}


\end{document}